\documentclass[11pt,twoside,a4paper]{article}
\usepackage{color,amscd,amsmath,amssymb,amsthm,latexsym,stmaryrd,authblk,mathabx,shuffle,hyperref,tikz,tkz-tab} 
\usepackage[latin1]{inputenc}
\usepackage[T1]{fontenc}   
\usepackage[english]{babel}
\usepackage[pdftex,all]{xy}
\usepackage{color}
\usepackage{esvect}
\usepackage{framed}
\usepackage{comment}
\usepackage{amsfonts}

\input{xy}
\xyoption{all}

\title{Lie and pre-Lie theory of Novikov algebras}
\date{}
\author{Ruggero Bandiera, Fr\'ed\'eric Patras}

\theoremstyle{plain}
\newtheorem{theo}{Theorem}[section]
\newtheorem{lemma}[theo]{Lemma}
\newtheorem{cor}[theo]{Corollary}
\newtheorem{prop}[theo]{Proposition}
\newtheorem{defi}[theo]{Definition}

\theoremstyle{remark}
\newtheorem{remark}{Remark}[section]
\newtheorem{notation}{Notation}[section]
\newtheorem{example}{Example}[section]

\newcommand{\K}{\mathbb{K}}

\newcommand{\sA}{\mathcal{A}}

\def\aaltbin#1#2{\ensuremath{\left(\kern-.35em\left(\genfrac{}{}{0pt}{}{#1}{#2}\right)\kern-.35em\right)} }

\renewcommand{\leq}{\leqslant}
\renewcommand{\geq}{\geqslant}
\newcommand{\N}{\mathbb{N}}
\newcommand{\Q}{\mathbb{Q}}
\newcommand{\R}{\mathbb{R}}

\newcommand{\tl}{\triangleleft}
\newcommand{\tr}{\triangleright}

\begin{document}
	
	\maketitle
	
	\begin{abstract}
		Novikov algebras provide a simple but powerful algebraic axiomatization of important features of classical differential calculus.
		We study their structure properties, modeling their relationships with commutative algebras with a derivation, featuring the role of their Lie and pre-Lie structures and analyzing the structure of their enveloping algebras. We focus on the combinatorial analysis of the Poincar\'e-Birkhoff-Witt (PBW) Theorem (classical and pre-Lie), the pre-Lie exponential and logarithm. The topic is important for applications of the theory and has been treated intensively for pre-Lie algebras. However, specific formulas can be obtained in the Novikov case. We analyze their structure, as well as featuring various remarkable properties. Related statistical phenomena on trees, tableaux and permutations are investigated in this context.
		
	\end{abstract}
	
	\tableofcontents
	
	\section{Introduction}
	
	Novikov algebras provide a simple but powerful algebraic axiomatization of important features of classical differential calculus \cite{DL2002}. They have attracted a renewed interest very recently due to their use in the context of stochastic differential equations and numerical methods. See e.g. \cite{bd,bkef,bh} and, for related developments in algebra \cite{diu,Manchon24,Foissy25}.
	
	In the present article
	we study their Lie-type structure properties, modeling their relationships with commutative algebras with a derivation, featuring the role of their Lie and pre-Lie structures and analyzing the structure of their enveloping algebras. The main focus is on the combinatorial analysis of the Poincar\'e-Birkhoff-Witt (PBW) theorem (classical and pre-Lie), the pre-Lie exponential, the pre-Lie logarithm and the associated flow map. The topic is important for applications of the theory and has been treated intensively for pre-Lie algebras. However, specific formulas can be obtained in the Novikov case. We analyze their structure, as well as featuring various remarkable properties. Statistical phenomena on trees, tableaux and permutations are investigated in this context.
	
	\vspace{0.5cm}
	
	The article is organised as follows. The first next three sections are mostly expositional and aim at giving a complete and self-contained account of the fundamental properties of Novikov algebras, especially those that can be obtained from their pre-Lie algebra structure.
	
	Section \ref{Novikov} surveys general definitions and properties related to Novikov algebras, free Novikov algebras, pre-Lie algebras and symmetric brace algebras.
	
	Section \ref{Diff} features various connections with algebras of differential operators, useful to understand how Novikov algebra structures can be used in the context of differential calculus.
	
	Section \ref{Env} details the structure of their enveloping algebras and the interpretation of computations in the enveloping algebra in terms of differential calculus.
	
	Section \ref{preLiePBW} investigates the consequences of the pre-Lie Poincar\'e-Birkhoff-Witt theorem.
	
	Section \ref{PBWc} studies the combinatorial and structural implications of the classical Poincar\'e-Birkhoff-Witt theorem.
	
	Section \ref{subsec:exp} investigates the pre-Lie exponential.
	
	Section \ref{sec:lie} investigates the pre-Lie logarithm.
	
	Section \ref{sec:flow} investigates the pre-Lie flow map.
	
	\vspace{0.5cm}
	
	We assume that the reader is familiar with the language and fundamental notions of the theory of pre-Lie algebras, coalgebras, bialgebras and Hopf algebras. Three references for the notions and results used throughout the article are \cite{Kapr,BS,Cartier2021}.\\
	
	{\bf Notation}\ \
	\begin{enumerate}
		\item Let $k\in \N$.  We denote by $[k]$ the set $\{1,\ldots,k\}$. By convention, $[0]=\emptyset$.  
		\item We work over the field of the rationals $\Q$. All results hold over an arbitrary field of characteristic $0$. 
		\item The $n$th symmetric group, that is the group of permutations of $[n]$,  is denoted $S_n$. We write $1_n$ its unit (the identity map). 
		\item The convolution product of linear endomorphisms of a bialgebra $B$, with product $m$, coproduct $\Delta$, unit $\eta$ and counit $\varepsilon$, is denoted $\ast$: 
		$$f\ast g:=m\circ (f\otimes g)\circ \Delta.$$
		It is associative, with unit denoted $\nu:=\eta\circ\varepsilon$. When the bialgebra is graded connected ($B=\bigoplus\limits_{n\in\N}B_n$, all structure maps preserve the graduation, and the degree 0 component identifies with the ground field), $\nu$ is the projection on the ground field orthogonally to the higher degree components. In that case, the bialgebra is automatically a Hopf algebra: its identity map has a convolution inverse $S$.
		\item We use the Sweedler notation and use $x^{(1)}\otimes x^{(2)}$ as a shortcut for $\Delta(x)$ for $x$ an arbitrary element in a coalgebra with coproduct $\Delta$. 
		\item A sequence $\mu=(\mu_1,\dots,\mu_k)$ of positive integers is called a composition of the integer $|\mu|:=\mu_1+\dots +\mu_k$ of lenght $l(\mu)=k$. 
		A weakly decreasing sequence is called a partition. We write 
		$$\mathcal P_n := \{(\mu_1,\dots , \mu_k)| \mu_1\geq \dots \geq \mu_k \geq 1, \mu_1 + \dots + \mu_k = n\}$$
the set of partitions 
of $n$. The set of partitions of $n$ in $k$ blocks is written $\mathcal P_{n,k}$.
		\item Partitions are equipped with the order defined by: $\mu=(\mu_1,\dots,\mu_k)< \nu=(\nu_1,\dots,\nu_l)$ if and only if $k\leq l$; for any $i\leq k$, $\mu_i\leq \nu_i$, and at least one of these $k+1$ inequalities is strict.
		\item The factorial $\mu!$ of a sequence $\mu=(\mu_1,\dots,\mu_k)$ of nonnegative integers is defined by $\mu!:=\prod\limits_{i=1}^k\mu_i!$.
		\item The cardinality of a set $S$ is denoted by $|S|$. For notational simplicity, when $f:S\to T$ is a map between sets and $t\in T$, we abbreviate $|f^{-1}(\{t\})|$ to $|f^{-1}(t)|$.
		\item We use as in \cite{Cartier2021} the term {\it gebra} for a vector space equipped with several (possibly many) algebraic structures. 
		\item Given a vector space $V$, we denote $\Q[V]$ the algebra of polynomials over $V$, that is, equivalently the algebra of covariants $\bigoplus\limits_{n\in \N}(V^{\otimes n})_{S_n}$. We use the monomial notation: given $v_1,\dots,v_n\in V$, $v_1\dots v_n$ stands for the class in the space of covariants $(V^{\otimes n})_{S_n}$ of the tensor $v_1\otimes\dots\otimes v_n$. Alternatively, when  the monomial notation could lead to ambiguities, we also denote the commutative monomial $v_1\dots v_n$ by  $v_1\odot\dots \odot v_n$.
		\item Given a tensor $v_1\otimes\dots\otimes v_n$ or a commutative monomial $v_1\dots v_n$ (depending of the context), given $I=\{i_1,\dots,i_k\}\subset [n]$ (where the $i_j$ are written in the natural order), we write $v_I$ for $v_{i_1}\otimes\dots\otimes v_{i_k}$, respectively $v_{i_1}\dots v_{i_k}$, with the convention $v_\emptyset:=1.$
	\end{enumerate}

	\section{Novikov algebras}\label{Novikov}
	
	A standard reference for fundamental properties of Novikov algebras is the article by Dzhumadil'daev and L\"ofwall \cite{DL2002}, to which we refer also for historical indications and the construction of free Novikov algebras, recalled below.
	
	\begin{defi} A right pre-Lie algebra (abbreviated to pre-Lie algebra in the present article) is a vector space $V$ equipped with a bilinear product $\triangleleft:V\otimes V\to V$ such that for all $x,y,z$ in $V$
		$$(x\tl y)\tl z-x\tl (y\tl z)=(x\tl z)\tl y-x\tl (z\tl y).$$ 
	\end{defi}
	
	The associated bracket $[x,y]_\tl:=x\tl y-y\tl x$ is a Lie bracket: it satisfies $[x,y]_\tl=-[y,x]_\tl$ and the Jacobi identity
	$$[[x,y]_\tl,z]_\tl+[[y,z]_\tl,x]_\tl +[[z,x]_\tl,y]_\tl=0.$$
	
	Left pre-Lie algebras are defined similarly: they are vector spaces with a bilinear product $\tr$ such that $(x\tr y)\tr z-x\tr (y\tr z)=(y\tr x)\tr z-y\tr (x\tr z)$. 
	The opposite algebra of a right pre-Lie algebra is left pre-Lie (setting $x\tr y:=y\tl x$). The associated Lie bracket $x\tr y-y\tr x=-[x,y]_\tl$ is denoted by $[x,y]_\tr$.
	
	Recall (see e.g. \cite{BS} and \cite[Section 6.2]{Cartier2021} for details and proofs) that the two notions of pre-Lie algebra and right symmetric brace algebra (abbreviated to symmetric brace algebra in the present article) are equivalent. 
	\begin{defi}
		A right symmetric brace algebra is a vector space $V$ equipped with a map 
		$$V\otimes \Q[V]\to V$$
		$$v\otimes P\longmapsto v\{P\}$$
		satisfying, for $v,w_1,\dots ,w_n,z_1,\dots ,z_m\in V$ the identities
		$$v\{1\}=v,$$
		\begin{equation}\label{sb1}
			(v\{w_1 \dots w_n\})\{z_1 \dots z_m\}=\sum\limits_{I_1\coprod\cdots\coprod I_{n+1}=[m]}v\{w_1\{ z_{I_1}\}\dots w_n\{ z_{I_n}\}z_{I_{n+1}}\},
		\end{equation}
		where some of the $I_i$ may be empty. 
	\end{defi}
	We will sometimes write  $v\{w_1, \dots ,w_n\}$ for $v\{w_1 \dots w_n\}$ to avoid notational ambiguities that might otherwise arise.
	
	For example,
	$$(v\{w_1w_2\})\{z\}=v\{w_1w_2z\}+v\{w_1\{z\}w_2\}+v\{w_1w_2\{z\}\}.$$
	The associated pre-Lie structure on $V$ is obtained by setting 
	\begin{equation}\label{sb2}v\tl w:=v\{w\}.\end{equation} 
	
	Conversely, the higher symmetric brace operations are recursively obtained from a pre-Lie product $\tl$, starting from the same relation, as 
	\begin{equation}\label{eqsybrace}
		v\{w_1 \dots w_n\}:=(v\{w_1 \dots w_{n-1}\})\{w_n\}-\sum\limits_{i=1}^{n-1}v\{w_1 \dots w_i\{w_n\}\dots w_{n-1}\}.
	\end{equation}
	
	Left symmetric brace algebras are defined similarly; the identity $\{w_1 \dots w_n\}v:=v\{w_1 \dots w_n\}$ puts in correspondence right and left symmetric brace algebra structures.
	
	\begin{defi}
		A right Novikov algebra (abbreviated Novikov algebra in the present article) is a pre-Lie algebra $(N, \tl)$ such that moreover
		$$a \tl (b  \tl  c) = b  \tl  (a  \tl  c)$$
		for all $a, b, c \in N$. \end{defi}
	
	Left Novikov algebras are defined similarly; the opposite algebra of a right Novikov algebra with product $\tl$ is a left Novikov algebra with product written $\tr$.
	
	\begin{example}[Commutative algebra with derivation]\label{comderex}
		If $\mathcal A$ is a commutative associative algebra and $\partial : \mathcal A \to \mathcal A$ is a derivation, the equation
		$a  \tl  b := \partial(a)b$ (resp. $a  \tr  b := a\partial(b)$)
		defines a Novikov (resp. left Novikov) algebra structure on $\mathcal A$. Indeed,
		$$a\tl (b\tl c)=b\tl (a\tl c)=\partial(a)\partial(b)c,$$ and
		$$(a\tl b)\tl c-a\tl (b\tl c)=\partial^2(a)bc+\partial(a)\partial(b)c-\partial(a)\partial(b)c=\partial^2(a)bc,$$ which is symmetric in $b$ and $c$.
		
		\item The corresponding symmetric brace algebra structure on $\mathcal A$ is then given by:
		\begin{equation}\label{identisymb}
			a\{a_1\dots a_n\}=a_1\dots a_n\partial^n(a)
		\end{equation}
		for $a,a_1,\dots,a_n\in \mathcal A$.
		We indeed have, by induction on $n$ and using Eq. (\ref{eqsybrace}):
		$$a\{a_1\dots a_n\}=(a\{a_1\dots a_{n-1}\})\{a_n\}-\sum\limits_{i=1}^{n-1}a\{a_1,\dots ,a_i\{a_n\},\dots ,a_{n-1}\}$$
		$$=(a_1\dots a_{n-1}\partial^{n-1}(a))\{a_n\}-\sum\limits_{i=1}^{n-1}a\{a_1,\dots ,a_n\partial(a_i) ,\dots ,a_{n-1}\}$$
		$$=a_n\partial(a_1\dots a_{n-1}\partial^{n-1}(a))-\sum\limits_{i=1}^{n-1}a_1\dots a_{i-1}a_n\partial(a_i)a_{i+1}\dots a_{n-1}\partial^{n-1}(a)$$
		$$=a_1\dots a_n\partial^n(a),$$
		by the Leibniz rule.
	\end{example}
	
	This example is generic, in the sense that free Novikov algebras can always be realized as sub Novikov algebras of commutative algebras with derivations:
	
	\begin{example}[Free Novikov algebras]\label{fnovex} 
		We have seen in Example \ref{comderex} that any commutative algebra with a derivation is a Novikov algebra. 
		Let now $A$ be a set, and $\mathcal A := \Q[a_i]$ the polynomial algebra generated by
		variables $a_i$ with $a\in{ A}, i\geq -1$.
		
		Given now $m = x_{i_1}y_{i_2}\dots  z_{i_k}$ a monomial in $\mathcal A$, we call
		$w(m) := i_1+\dots+i_k$ its weight and $d(m) := k$ its degree. For example, $w(a_{-1}^3a_{5}b_{0}^4)=-3+5+0=2$ and $d(a_{-1}^3a_{5}b_{0}^4)=3+1+4=8$.
		
		Finally, define a derivation $\partial : \mathcal A \to \mathcal A$
		according to $\partial (a_i) := a_{i+1}$, and let  $\tl$ be the corresponding Novikov product on $\mathcal A$, defined as in
		the previous example. We denote by $\mathcal N(A) \subset \mathcal A$ the subspace spanned by monomials
		of weight -1. One easily checks that $(\mathcal N(A),  \tl )$ is a Novikov subalgebra of $(\mathcal A,  \tl )$. In fact,
		it can be proved that it is the free Novikov algebra generated by $A$ (\cite{DL2002} Thm 7.8). 
		
	\end{example}
	\begin{example}[Free Novikov algebra on one generator]\label{fnovex1gen}
		In the particular case $A = \{x\}$, we denote by $\mathcal N (x)$ the free Novikov algebra generated by $x$. It is graded
		by polynomial degree, and we denote by $\mathcal N_n (x)$  the component of degree $n$. 
		Given $\mu = (\mu_1,\dots  , \mu_k) \in \mathcal P_n$, we denote by $x_\mu \in \mathcal N_{n+1}(x)$ the
		monomial
		$$x_\mu := x_{\mu_1-1}\dots x_{\mu_k-1}x_{-1}^{n-k+1}.$$
		For $n = 0$ we denote by $\mathcal P_0$ the set containing only the empty partition, which we denote by
		$(0)$, and we set $x_{(0)} := x_{-1}$. When $\mu$ varies in $\mathcal P_n$, the monomials $x_\mu$ form a vector space basis
		of $\mathcal N_{n+1}(x)$.
	\end{example}
	
	Let us conclude with a link to free pre-Lie algebras. The free pre-Lie algebra over a generator $x$, $p\mathcal{L}(x)$, is also the free symmetric brace algebra over $x$ and has therefore a graded ordered basis $B=\coprod\limits_{n\in\N^\ast}B_n$:
	$B_1=\{x\}$, $B_2=\{x\{x\}\}$, $B_3=\{x\{x,x\},x\{x\{x\}\}\}\dots$.
	In general, $B_{n+1}$ is recursively obtained from the $B_i,\ i\leq n$ as follows. Assume that an order has been defined on each $B_i$ and that $B_j<B_k$ when $j<k$. 
	Then, write $|b|=k$ if $b\in B_k$. Write $X_n$ for the lexicographically ordered set of all weakly decreasing sequences $(z_1,\dots ,z_k)$ of elements of $\coprod\limits_{i\leq n}B_i$ such that
	$|z_1|+\dots + |z_k|=n$. Define then
	$$B_{n+1}:=\{x\{z_1\dots z_k\}, (z_1,\dots ,z_k)\in X_n\}.$$
	Use the lexicographical order on $X_n$ to define a total order on $B_{n+1}$ (by requiring $x\{y\}\leq x\{z\}$ for $y,z\in X_n$ if and only if $y\leq z$); finally define a total order on $\coprod\limits_{k\leq n+1}B_k$ by requiring that $B_j<B_{n+1}$ for $j\leq n$.
	This construction is a variant of  the ones in \cite{AG,BS,Cartier2021}, to which we refer for details and a proof that the basis so obtained is indeed a basis of the free pre-Lie algebra or, equivalently, of the free symmetric brace algebra. 
	
	Equivalently, the free pre-Lie algebra is parametrized by (non planar) rooted trees \cite{Ermolaev,Chapoton2001}. The bijection $T$ between $B$ and the set of trees is recursively obtained  by setting $T(x):=\bullet$ and mapping $x\{z_1\dots z_k\}$ to the tree with a root $\bullet$ to which are attached the trees $T(z_1),\dots, T(z_k)$. For example, $x\{x,x\{x,x\{x\}\}\}$ is mapped to the tree
	
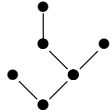
\[  \xy 
	{\ar@{-}(0,-4)*{\bullet};(-4,0)*{\bullet}};
	{\ar@{-}(0,-4)*{\bullet};(4,0)*{\bullet}};
	{\ar@{-}(4,0)*{\bullet};(0,4)*{\bullet}};
	{\ar@{-}(4,0)*{\bullet};(8,4)*{\bullet}};
	{\ar@{-}(0,4)*{\bullet};(0,9)*{\bullet}};
	
	\endxy\]
	
	As a Novikov algebra is pre-Lie, the injection of $\{x\}$ into ${\mathcal N}(x)$, the free Novikov algebra generated by $x$, induces a (canonical) map $\rho$ from $p\mathcal{L}(x)$, the free pre-Lie algebra over $x$, to ${\mathcal N}(x)$. Using the fact that a symmetric brace operation in the free pre-Lie algebra is mapped to the same symmetric brace operation interpreted in the free Novikov algebra, it follows from Eq (\ref{identisymb}) and the construction of ${\mathcal N}(x)$ in Example \ref{fnovex1gen} 
	that this map is recursively obtained, using our previous notation, as:
	$$\rho(x)=x_{-1},$$
	$$\rho(x\{z_1\dots z_k\})=x_{k-1}\rho(z_1)\dots \rho(z_k).$$
	Similarly, a tree with $n$ vertices, viewed as an element of $p\mathcal{L}(x)$, is mapped to the product $\prod\limits_{v}x_{f(v)-1}$, where $v$ runs over the vertices of the tree (including the root and the leaves) and $f(v)$ stands for the number of outgoing upper branches (in some references this is called the \emph{fertility} of the vertex $v$).
	For example, $$\rho(x\{x\{x,x,x\},x\{x\}\})=x_1x_2x_{-1}^3x_0x_{-1}=x_0x_1x_2x_{-1}^4.$$
	
	\section{Relations to differential calculus}\label{Diff}\label{comderalg}
	
	We illustrate now an important point: the deep relationships between Novikov algebras and differential operators.
	
	Let $(\mathcal A,\partial)$ be a commutative algebra with a derivation, as in Example \ref{comderex} above. Write $FDiff(\mathcal A)$ and call {\it gebra of formal differential operators on $(\mathcal A,\partial)$} the vector space $\mathcal A[D]=\bigoplus_{n\in\N}\mathcal A\cdot D^n$, where $D^0:=1$ and $D$ stands for a dummy variable. 
	The gebra $FDiff(\mathcal A)$ carries the following structures:
	\begin{enumerate}
		\item It is a graded vector space and graded left $\mathcal A$-module, with degree $n$ component $FDiff_n(\mathcal A)$ spanned (as an $\mathcal A$-module) by $D^n$,
		\item it is a commutative graded algebra: the $\mathcal A$-algebra of polynomials in $D$, with product denoted $\bullet$,
		$$(a\cdot D^n)\bullet (b\cdot D^m):=ab\cdot D^{n+m},$$
		\item it is an associative algebra for the product $\ast$ inductively defined by
		$$(a\cdot D^n)\ast (b\cdot D^m):=(a\cdot D^{n-1})\ast (\partial(b)\cdot D^m+b\cdot D^{m+1}),$$
		where $n\geq 1$ and the convention $a\ast (b\cdot D^m):=ab\cdot D^m$, or, equivalently, by
		\begin{equation}\label{leib1}
			(a\cdot D^n)\ast (b\cdot D^m):=\sum\limits_{p=0}^n{n\choose p}a\partial^{n-p}(b)D^{m+p}.
		\end{equation}
		Associativity follows from Leibniz rules for derivations. Notice that $(a\cdot D^n)\bullet (b\cdot D^m)=(a\cdot D^n)\ast (b\cdot D^m)+l.o.t,$
		where l.o.t. (lower order term) stands for an element of total degree less than $n+m$.
		\item The degree 1 component $FDiff_1(\mathcal A)$ is a left Novikov algebra for the product $a\cdot D\tr b\cdot D:=a\partial(b)\cdot D$. It is obviously isomorphic via $$a\cdot D\longmapsto a$$ to the left Novikov algebra $(\mathcal A,\tr)$ of Example \ref{comderex}. 
		\item The associated Lie bracket is obtained as
		$$[a\cdot D,b\cdot D]_\tr:=(a\partial(b)-b\partial(a))\cdot D,$$
		it makes $FDiff_1(\mathcal A)$ a sub Lie algebra of $FDiff(\mathcal A)$ when the latter is equipped with the Lie bracket obtained by antisymmetrization of the associative product $\ast$,
		$$[a\cdot D^n,b\cdot D^m]_\ast:=(a\cdot D^n)\ast (b\cdot D^m)-(b\cdot D^m)\ast (a\cdot D^n).$$
		Indeed, $(a\cdot D)\ast (b\cdot D)-(b\cdot D)\ast (a\cdot D)=a\partial(b)\cdot D+ab\cdot D^2-b\partial(a)\cdot D-ba\cdot D^2,$
		and thus, on $FDiff_1(\mathcal A)$ we have $$[a\cdot D,b\cdot D]_\tr=[a\cdot D,b\cdot D]_\ast.$$
		
	\end{enumerate}
	
	\begin{remark} Recall Grothendieck's recursive construction of the subalgebra $Diff(\sA)\subset End(\sA)$ of differential operators on $\sA$. Given $a\in\sA$, we denote by $\mu_a:\sA\to\sA:x\to ax$ the operator of multiplication by $a$. One defines $Diff_{\le-1}(\sA)=0$, and for $k\ge0$ one defines the subspace $Diff_{\le k}(\sA)\subset End(\sA)$ of \emph{differential operators on $\sA$ of order $\le k$} recursively via
		\[ Diff_{\le k}(\sA) = \{ D\in End(\sA)\,\,\mbox{s.t.}\,\, [D,\mu_a]\in Diff_{\le k-1}(\sA),\,\forall\,a\in\sA\}.\] 
		For instance, $\phi\in Diff_{\le0}(\sA)$ if and only if $\phi(ab)-a\phi(b)=0$ for all $a,b\in\sA$: taking $b=\mathbf{1}_\sA$ we get $\phi(a)=a\phi(\mathbf{1}_\sA)$, thus $\phi=\mu_{\phi(\mathbf{1}_\sA)}$, showing that $Diff_{\le0}(\sA)$ is the subalgebra of $End(\sA)$ generated by the multiplication operators $\mu_a$.
		
		Similarly, $\phi\in Diff_{\le1}(\sA)$ if for all $a\in\sA$ there exists a $c\in\sA$ (depending on $a$) such that $\phi(ab)-a\phi(b)=cb$ for all $b\in\sA$: taking $b=\mathbf{1}_\sA$ we get $c=\phi(a)-a\phi(\mathbf{1}_\sA)$, and thus $\phi(ab)-a\phi(b)-b\phi(a)+ab\phi(\mathbf{1}_{\sA})$ for all $a,b\in\sA$. In particular, derivations of $\sA$ can be characterized as the subspace $Der(\sA)\subset Diff_{\le1}(\sA)$ of differential operators of order $\le1$ vanishing on the unit $\mathbf{1}_\sA$. 
		
		The Poisson and Jacobi identities \[[D\circ D',\mu_a]=D\circ[D',\mu_a]+[D,\mu_a]\circ D',\qquad[[D,D'],\mu_a]=[D,[D',\mu_a]]+[[D,\mu_a],D'],\] together with an obvious induction on $h,k$, show that  $Diff_{\le h}(\sA)\circ Diff_{\le k}(\sA)\subset Diff_{\le h+k}(\sA)$ and $[Diff_{\le h}(\sA),Diff_{\le k}(\sA)]\subset Diff_{\le h+k-1}(\sA)$ for all $h,k\ge0$. This shows that $Diff(\sA)=\bigcup_{k\ge0}Diff_{\le k}(\sA)\subset End(\sA)$ is a subalgebra with respect to the composition product, and is a filtered algebra with the natural filtration $Diff_{\le 0}(\sA)\subset\cdots\subset Diff_{\le k}(\sA)\subset\cdots\subset Diff(\sA)$: in particular $a\partial^n\in Diff_{\le n}(\sA)$ for all $a\in\sA$, $n\ge0$. There is thus a natural map \[ev:FDiff(\sA)\to Diff(A):aD^n\to a\partial^n,\]
		which is easily checked to be an algebra map, since the product $\ast$ on $FDiff(\sA)$ was implicitly defined using the Leibniz rule (Eq (\ref{leib1})). In particular, this restricts to $\sA=FDiff_0(\sA)\to Diff_{\le0}(\sA)\subset Diff(\sA):a\to \mu_a$, making $Diff(\sA)$ an $\sA$-algebra, and the image of $ev$ can be characterized as the $\sA$-subalgebra of $Diff(\sA)$ generated by the derivation $\partial$. When $a\to\mu_a$ is injective and the powers $\partial^n$ are linearly independent over $\sA$ the map $ev$ is also injective, and it can be used to identify $FDiff(\sA)$ with the aforementioned $\sA$-subalgebra of $Diff(\sA)$.
	\end{remark}
	
	\begin{example} Particularly interesting examples of the previous construction are when $(\sA,\partial)$ is the algebra $(C^\infty(\R),d)$ of smooth functions on the real line, or the algebra of polynomials $(\R[X],d)$, where $d$ stands for the usual derivation operator. In these cases $FDiff(\sA)$ identifies with the algebra of (smooth or polynomial, respectively) differential operators on the real line, while $FDiff_1(\sA)$ identifies with the Novikov (and pre-Lie, and Lie) algebra of (smooth or polynomial, respectively) vector fields.
	\end{example}

	\section{Enveloping algebras of Novikov algebras}\label{Env}
	We make now explicit the structure of enveloping algebras of Novikov algebras. We choose $\Q$ as a ground field, as usual.
	
	Recall, for completeness sake, that if $L$ is a Lie algebra, an enveloping algebra of $L$ is an associative algebra, usually denoted $U(L)$ such that there is a natural isomorphism $Hom_{Lie}(L,B)\cong Hom_{Ass}(U(L),B)$, where $B$ runs over associative algebras and is equipped with the Lie algebra structure $[b,b']=bb'-b'b$, and where $Hom_{Lie}$ and $Hom_{Ass}$ stand respectively for the set of Lie algebra and associative algebra morphisms. The morphism $l\longmapsto l\otimes 1+1\otimes l$ induces a map $\Delta$ from $U(L)$ to $U(L)\otimes U(L)$ which is an algebra map and equips $U(L)$ with a Hopf algebra structure. The Lie algebra $L$  identifies then to the sub Lie algebra of primitive elements of $U(L)$, $Prim(U(L)):=\{x\in U(L),\Delta(x)=x\otimes 1+1\otimes x\}$. By the Poincar\'e-Birkhoff-Witt (PBW) Theorem, $(U(L),\Delta)$ is isomorphic to $(\Q[L],\Delta)$ as a (conilpotent cocommutative cofree) coalgebra, where $\Q[L]$ is equipped with the unshuffle coproduct:
	$$\Delta(l_1\dots l_n):=\sum\limits_{I\coprod J=[n]}l_I\otimes l_J.$$ 
	Enveloping algebras are defined up to isomorphism. The most standard construction, that we will denote $U_{st}(L)$, is obtained by considering the quotient of the tensor algebra $T(L):=\bigoplus\limits_{n\in\N}L^{\otimes n}$ by the ideal generated by elements of the form $l\otimes l'-l'\otimes l-[l,l'],\ l,l'\in L$.
	Further details on enveloping algebras and the PBW theorem can be found in \cite{Reutenauer,Cartier2021} or any other reference book on the subject.

	Another construction of the enveloping algebra of a Novikov algebra $N$ is easily obtained from its symmetric brace algebra structure. This is a particular case of the construction of the enveloping algebra of a pre-Lie algebra \cite[Chap. 6]{Cartier2021}.
	Let us set: $(U(N),\Delta):=(\Q[N],\Delta)$ as a coalgebra and define the product $\ast$ making it a Hopf algebra and the enveloping algebra of $N$ by:
	$\forall w_1,\dots,w_n,z,z_1,\dots, z_m\in N$,
	$$1\ast z=z\ast 1=z,$$
	and
	\begin{equation}\label{envprod}
		w_1\dots w_n \ast z_1\dots z_m:=\sum\limits_{I_1\coprod\dots\coprod I_{n+1}=[m]}w_1\{z_{I_1}\}\dots  w_n\{z_{I_n}\}z_{I_{n+1}},
	\end{equation}
	where some of the $I_i$ may be empty.
	Notice that this implies:
	\begin{equation}\label{envprod2}
		w_1\dots w_n \ast z_1\dots z_m:=\sum\limits_{I_1\coprod I_2}w_1\{z_{I_1}\}(w_2\dots w_n \ast z_{I_2}).
	\end{equation}
	
	\begin{example}
		For example, we have:
		$$w_1w_2\ast z_1z_2=(w_1\{z_1z_2\})w_2+w_1(w_2\{z_1z_2\})+(w_1\{z_1\})(w_2\{z_2\})+(w_1\{z_2\})(w_2\{z_1\})$$
		$$+(w_1\{z_2\})w_2z_1+w_1z_1(w_2\{z_2\})+(w_1\{z_1\})w_2z_2+w_1(w_2\{z_1\})z_2+w_1w_2z_1z_2,$$
		where we have put parentheses to clearly identify the elements of $N$ appearing in the monomials on the right side of the equation.
	\end{example}
	
	\begin{lemma}\label{compdelta}
		We have, for $w,z\in U(N)$,
		$$\Delta(w\ast z)=\Delta(w)\ast\Delta(z).$$ 
	\end{lemma}
	
	\begin{proof}
		Indeed, let us assume, without restriction, that $w=w_1\dots w_n$ and $z=z_1\dots z_m$, with our previous notation. Then,
		$$\Delta(w\ast z)=\sum\limits_{A\coprod B=[n]\atop I_1\coprod\dots\coprod I_{n+1}=[m] }\sum\limits_{C\coprod D=I_{n+1}}w_{a_1}\{z_{I_{a_1}}\}\dots w_{a_{|A|}}\{z_{I_{a_{|A|}}}\}z_C\otimes w_{b_1}\{z_{I_{b_1}}\}\dots w_{b_{|B|}}\{z_{I_{b_{|B|}}}\}z_D,$$
		with the notation $A=\{a_1,\dots,a_{|A|}\},\ B=\{b_1,\dots,b_{|B|}\}.$ Grouping together terms such that $I_{a_1}\coprod\dots \coprod I_{a_{|A|}}\coprod C=U$ and  $I_{b_1}\coprod\dots \coprod I_{b_{|B|}}\coprod D=V$ with $U\coprod V=[m]$, one gets that $\Delta(w\ast z)$ is equal to
		$$\sum\limits_{A\coprod B=[n]\atop U\coprod V=[m]}(w_A\ast z_U)\otimes (w_B\ast z_V)=\left( \sum\limits_{A\coprod B=[n]}w_A\otimes w_B\right)\ast \left( \sum\limits_{U\coprod V=[m]}z_U\otimes z_V\right)$$
		$$=\Delta(w)\ast\Delta(z).$$
	\end{proof}
	
	\begin{prop}
		The triple $(U(N),\ast,\Delta)$ together with the canonical embedding of $N$ into $\Q[N]=U(N)$ is an enveloping algebra of $N$.
	\end{prop}
	
	\begin{proof}
		Assume first that the product $\ast$ is associative. Notice that for $x,y\in N\subset\Q[N]$ we have $x\ast y -y\ast x = xy +x\{y\}-yx-y\{x\} = x\{y\}-y\{x\}=[x,y]_{\triangleleft}$. Then $N\subset\Q[N]$ is a Lie subalgebra, and by the universal property of enveloping algebras, if $A$ is an enveloping algebra for $N$ the identity map from $N\subset A$ to $N\subset\Q[N]$ induces an associative algebra map $\phi$ from $A$ to $(\Q[N],\ast)$. Consider now the diagram 
		\[ \xymatrix{ A \ar[r]^-\phi\ar[d]_-{\Delta} & \Q[N]\ar[d]^-{\Delta} \\ A\otimes A\ar[r]_-{\phi\otimes\phi} & \Q[N]\otimes\Q[N], } \] 
		where all maps are algebra morphisms. Since $N$ sits inside both $A$ and $\Q[N]$ as the subspace of primitive elements, the restrictions of $\Delta\circ\phi$ and $(\phi\otimes\phi)\circ\Delta$ to $N$ coincide: hence, the universal property of $A$ implies that the diagram is commutative, so $\phi$ is also a coalgebra morphism. As, by the PBW theorem, $A$ is isomorphic to $\Q[N]$ as a (conilpotent cocommutative cofree) coalgebra, and since $\phi$ restricts to the identity on $N$, this map $\phi$ is necessarily a coalgebra isomorphism and thus also an algebra and Hopf algebra isomorphism. Here we used the fact that a coalgebra morphism $F:\mathbb{Q}[V]\to\mathbb{Q}[W]$ is an isomorphism if and only if so is the map of vector spaces $f:V\hookrightarrow\mathbb{Q}[V]\xrightarrow{F}\mathbb{Q}[W]\twoheadrightarrow W$ (the first and last arrow in the composition being the natural inclusion and projection respectively): this can be seen as a coalgebraic analogue of the usual inverse function Theorem.
		
		The proof of the Proposition will thus follow if we show that the product $\ast$ is associative. 
		Let us consider $w=w_1\dots w_n,\ z=z_1\dots z_m,\ y=y_1\dots y_k$ with the $w_i,z_i,y_i$ in $N$. We argue by induction on triples $(n,m,k)$ equipped with the obvious (partial) order. The case when $n$ or $m$ or $k=0$ corresponds to assuming $w=1$ or $z=1$ or $y=1$ and is straightforward. We then have, using Eq (\ref{envprod2}), the Sweedler notation and Lemma \ref{compdelta}:
		$$(w\ast z)\ast y=\Big(w_1\{z^{(1)}\}(w_2\dots w_n\ast z^{(2)})\Big)\ast y$$
		$$=(w_1\{z^{(1)}\})\{y^{(1)}\}\Big(w_2\dots w_n\ast z^{(2)}\ast y^{(2)}\Big)$$
		$$=w_1\{(z\ast y)^{(1)}\}(w_2\dots w_n\ast (z\ast y)^{(2)}),$$
		$$=w\ast(z\ast y).$$
	\end{proof}
	
	\begin{remark} The enveloping algebra of a left Novikov algebra is defined similarly, with an associative product given by the opposite product $\ast^{op}$:
		\begin{equation}\label{envprodopp}
			z_1\dots z_m\ast^{op}  w_1\dots w_n :=\sum\limits_{I_1\coprod\dots\coprod I_{n+1}=[m]}\{z_{I_1}\}w_1\dots  \{z_{I_n}\}w_nz_{I_{n+1}},
		\end{equation}
		All properties and constructions related to enveloping algebras dualize in an obvious way from Novikov algebras to left Novikov algebras.
	\end{remark}
	
	Let us apply now these constructions to differential operators. Let $(\mathcal A,\partial)$ be a commutative algebra with a derivation. We view it as a left Novikov algebra. The map $a\longmapsto a\cdot D$ from $\mathcal A$ to $FDiff_1(\mathcal A)$ is a morphism of left Novikov and of Lie algebras and induces therefore an algebra map $\gamma$ from $(U(\mathcal A),\ast^{op})$ to $(FDiff,\ast)$. We write abusively $ev:U(\sA)\to Diff(\sA)$ for the map induced by left composition of $\gamma$ with the  map $ev:FDiff(\sA)\to Diff(\sA)$ defined in the previous section.
	
	\begin{prop}\label{mappingrule}
		The map $\gamma$, which is an algebra map from $(U(\mathcal A)=\Q[\mathcal A],\ast^{op})$ to $(FDiff,\ast)$ is also an algebra map from $(U(\mathcal A)=\Q[\mathcal A],\odot)$ to $(FDiff,\bullet)$, where, for notational clarity, $ \odot$ denotes here the product of polynomials. In particular, for $a_1 \odot \dots  \odot a_n$ a monomial in $\Q[\mathcal A]$
		$$\gamma(a_1 \odot \dots  \odot a_n)=a_1\dots a_n\cdot D^n,$$
		and therefore
		$$ev(a_1 \odot \dots  \odot a_n)=a_1\dots a_n\partial^n.$$
	\end{prop}
	\begin{proof}
		Let us prove the Proposition by induction on $n$. We have:
		$$\gamma(a_1 \odot \dots  \odot a_n)=\gamma(a_1\ast^{op}(a_2\odot \dots  \odot a_n)-\sum\limits_{i=2}^n a_2 \odot \dots  \odot \{a_1\}a_i \odot\dots  \odot a_n)$$
		$$=a_1D\ast \gamma(a_2\odot \dots \odot a_n)-\sum\limits_{i=2}^n\gamma(a_2 \odot \dots  \odot a_1\partial a_i \odot\dots  \odot a_n)$$
		$$=a_1D\ast (a_2\dots  a_nD^{n-1})-\sum\limits_{i=2}^n a_1a_2\dots \partial a_i\dots  a_n D^{n-1}$$
		$$=a_1a_2\dots  a_nD^{n}+\sum\limits_{i=2}^n a_1a_2\dots \partial a_i\dots  a_n D^{n-1}-\sum\limits_{i=2}^n a_1a_2\dots \partial a_i\dots  a_n D^{n-1}$$
		$$=a_1a_2\dots  a_nD^{n}.$$
	\end{proof}
	
	\section{Pre-Lie PBW Theorem in Novikov algebras}\label{preLiePBW}
	
	Given $(L,\triangleleft)$ a pre-Lie algebra, consider its enveloping Hopf algebra $(\Q[L],\ast,\Delta)$ constructed as for Novikov algebras in the previous section, and a second enveloping Hopf algebra of $L$, $(U(L),\cdot,\Delta)$, for instance the standard one, $U_{st}(L)$, whose elements are classes of tensors. 
	
	There is then an isomorphism of Hopf algebras $\operatorname{pbw}_\triangleleft^{-1}:U_{st}(L)\to\Q[L]$ defined by 
	$$\operatorname{pbw}_\triangleleft^{-1}(x_1\otimes\cdots \otimes x_n)=x_1\ast\cdots\ast x_n$$ for $x_1,\ldots,x_n\in L.$
	We denote the inverse by $\operatorname{pbw}_\triangleleft:\Q[L]\to U(L)$, and call it the \emph{pre-Lie Poincar\'e-Birkhoff-Witt isomorphism}: it is an isomorphism of Hopf algebras, and in particular of coalgebras, but it {\it does not} coincide with the usual Poincar\'e-Birkhoff-Witt isomorphism, to be studied in section \ref{PBWc}. 
	
	We denote by $\eta:U(L)\to L$ the corestriction of $\operatorname{pbw}_\triangleleft^{-1}$, that is, its composition with the natural projection $p:\Q[L]\to L$. Since $\operatorname{pbw}_\triangleleft^{-1}$ is a morphism of coalgebras and $\Q[L]$ is cofree (as a conilpotent cocommutative coalgebra), we can reconstruct $\operatorname{pbw}_\triangleleft^{-1}$ from $\eta$ according to 
	(see for example \cite[Section 2.13]{Cartier2021})
	\begin{equation}\label{pdfFlaPBW}
		\operatorname{pbw}_\triangleleft^{-1}(X) = \sum_{k\ge1}\frac{1}{k!}\eta(X^{(1)})\odot\cdots\odot\eta(X^{(k)}), 
	\end{equation}
	where we use Sweedler's notation $\Delta^{k-1}(X)=X^{(1)}\otimes\cdots\otimes X^{(k)}$ for the iterated coproduct $\Delta^{k-1}:U(L)\to U(L)^{\otimes k}$. The definition of the product 
	$\ast$ on $\Q[L]$ implies that $p(Y\ast x) =p(Y)\{x\}=p(Y)\triangleleft x$ for all $x\in L$ and $Y\in\Q[L]$, and thus
	\begin{equation}\label{productFlaPBW}
		\eta(x_1\otimes\cdots \otimes x_n) = p\left(x_1\ast\cdots\ast x_n\right)= (\cdots(x_{1}\triangleleft x_{2})\cdots)\triangleleft x_{n}.
	\end{equation}
	Given a totally ordered set $I=\{i_1<\cdots< i_h\}$, we shall use the shorthand
	\[ x^{\triangleleft I}:=(\cdots(x_{i_1}\triangleleft x_{i_2})\cdots)\triangleleft x_{i_h}.\]
	With this notation, formula \eqref{pdfFlaPBW} becomes 
	\[ \operatorname{pbw}_\triangleleft^{-1}(x_1\otimes\cdots \otimes x_n) = x_1\ast\cdots\ast x_n=\sum_{k=1}^n\sum_{I_1\bigsqcup_< \cdots \bigsqcup_< I_k=[n]} x^{\triangleleft I_1}\odot\cdots\odot x^{\triangleleft I_k},\]
	where $I_1\bigsqcup_< \cdots \bigsqcup_< I_k=[n]$ means that the blocks of the partition are ordered according to $\min(I_1)<\cdots<\min(I_k)$.
	
	Similar considerations can be repeated for a left pre-Lie algebra $(L,\triangleright)$, in which case we obtain the isomorphism $\operatorname{pbw}_\triangleright^{-1}:U(L)\to \Q[L]$ 
	\begin{equation}\label{flapbwpl}
		\operatorname{pbw}_\triangleright^{-1}(x_1\otimes\cdots \otimes x_n) = x_1\ast^{op}\cdots\ast^{op} x_n=\sum_{k=1}^n\sum_{I_1\bigsqcup_< \cdots \bigsqcup_< I_k=[n]} x^{\triangleright I_1}\odot\cdots\odot x^{\triangleright I_k},
	\end{equation}
	where $x^{\triangleright I}:= x_{i_1}\triangleright(\cdots(x_{i_{h-1}}\triangleright x_{i_h})\cdots)$ if $I=\{i_i<\cdots<i_h\}$.
	\begin{example}\label{pLtoNov}
		Consider the generic case when $L$ is a free left  pre-Lie algebra on $n$ generators $x_1,\ldots,x_n$: it has a vector space basis indexed by isomorphism classes of pairs $(T,\ell)$, where $T$ is a rooted tree and $\ell:V(T)\to[n]$ is a labeling of the vertices of $T$  by elements of $[n]$, and where each such labeled tree is regarded as an iterated brace $b_{(T,\ell)}\in L$ in the generators, as explained in section \ref{Novikov}. Denote by $\mathcal{T}^{dec}_n$ the set of pairs $(T,\ell)$ as above, such that moreover $\ell:V(T)\to[n]$ is bijective and decreasing, meaning that the label of a vertex is always greater than the labels of its descendants. In this case  $\eta(x_1\otimes\cdots \otimes x_n)=x^{\triangleright[n]}$ can be expanded as $\eta(x_1\otimes\cdots \otimes x_n)=\sum_{(T,\ell)\in\mathcal{T}^{dec}_n}b_{(T,\ell)}$, and similarly, by (\ref{flapbwpl}), $\operatorname{pbw}_\triangleright^{-1}(x_1\otimes\cdots \otimes x_n)$. can be expanded as a sum of decreasing forests with $n$ vertices bijectively labeled by $[n]$. 
	\end{example}
	
	The aim of the remaining part of this section is to explicitly describe $\eta:U(A)\to A$ when $(A,\triangleright)$ is the left Novikov algebra associated with a commutative algebra $A$ with a derivation $\partial\in\operatorname{Der}(A)$, in which case we have
	\[ \eta(f_1\otimes\cdots \otimes f_i) =f^{\triangleright[i]}= (f_1\partial)\circ\cdots\circ(f_{i-1}\partial)(f_i). \]
	As we have seen earlier, as free Novikov algebras can be realized as sub Novikov algebras of a commutative algebra $A$ with a derivation, this case is generic. We first get from the general formula for pre-Lie algebras:
	\begin{lemma}\label{lempbwpl}
		We have
		\[ \eta(f_1\otimes\cdots \otimes f_i) = \sum_{(T,\ell)\in\mathcal{T}_i^{dec}}\prod_{v\in V(T)} \partial^{\operatorname{val}(v)}(f_{\ell(v)}), \]
		where $\operatorname{val}(v)$ is the number of children of $v\in V(T)$. 
	\end{lemma}
	
	\begin{notation} We shall denote by $Inc(i)$ the set of maps $\phi$ from $[i-1]$ to $[i]$ such that $\phi(k)>k$ for all $k\in[i-1]$ (notice in particular that we always have $|\phi^{-1}(1)|=0$).	\end{notation}
	
The above formula can be obtained more directly as follows.

	Using the Leibniz formula, one gets by induction on $i$ that
	\begin{equation}\label{eq:eta} (f_1\partial)\circ\cdots\circ(f_{i-1}\partial)(f_i)=\sum\limits_{\phi\in Inc(i)}f_1^{(|\phi^{-1}(1)|)}f_2^{(|\phi^{-1}(2)|)}\dots f_i^{(|\phi^{-1}(i)|)},
	\end{equation}
	Here, $\phi(k)=l$ means that the derivative following $f_k$ in the expression ``hits'' $f_l$ when the Leibniz formula is applied to get 
	$$f_k\partial ((f_{k+1}\partial)\circ\cdots\circ(f_{i-1}\partial)(f_i))$$
	$$=f_k(f'_{k+1}\partial)\circ\cdots\circ(f_{i-1}\partial)(f_i)+\dots +f_k(f_{k+1}\partial)\circ\cdots\circ(f_{i-1}\partial)(f'_i).$$
	The labeled tree in $\mathcal{T}_i^{dec}$ bijectively associated to $\phi$ is then the tree with edges $(i,\phi(i))$.
	
	
	\ \par
	
	Now, given $\phi\in Inc(i)$, it follows from the definitions that  $|\phi^{-1}(1)|+\dots +|\phi^{-1}(k)|\leq k-1$ for all $k\in[i-1]$. This motivates the following definition.
	
	\begin{notation} Given $i\ge1$, we shall denote by $K_i$ the set of length $i$ words $w=w_1\cdots w_i$ in the alphabet $\mathbb{N}$ such that $w_1+\cdots+w_j\le j-1$ for all $1\le j<i$ and $w_1+\cdots +w_i= i-1$. Notice that in particular $w_1=0$.
	\end{notation}	
	For instance
	\[ K_1 = \{0\},\qquad K_2 = \{01\},\qquad K_3 = \{002,011\},\qquad K_4 = \{ 0003, 0012, 0021, 0102, 0111\},\]
	\[ K_5 = \{ 00004, 00013, 00022, 00031, 00103, 00112, 00121, 00202\}\]
	\[\coprod\{00211, 01003,01012, 01021, 01102,01111\}\]
	
	Up to dropping the first zero and reverting the order of letters in words, this family of sets identifies with exemple (86) of Catalan objects in Stanley's book \cite{Stanley}: sequences $(a_1, \dots ,a_n)$ of nonnegative integers satisfying $a_1+\cdots +a_i \geq i$ with $a_1+\cdots +a_n=n$.
	In particular, the cardinalities of these sets are the Catalan numbers $|K_{i+1}|= \frac{1}{i+1}\binom{2i}{i}$: this will also be seen directly in the following Lemma \ref{freemagma}, by establishing an explicit bijective correspondence between $K_{i}$ and the set of planar binary rooted trees with $i$ leaves.
	
	Given $w=w_1\cdots w_i\in K_i$, it will  be convenient to introduce the auxiliary word $d(w)=d_1\cdots d_i$, where for all $1\le j\le i$ we set \[d_j:=j-1-w_1-\cdots-w_j.\]
	We then associate a coefficient $c_w\in\mathbb{N}$ to $w$ according to the formula 
	\[ c_w :=\prod_{j=1}^i\binom{d_j+w_j}{w_j} \]
	(notice that $\binom{d_j+w_j}{w_j}=1$ whenever $w_j=0$ or $d_j=0$, so many of the factors don't actually contribute to the product). For instance, for $w=0002013\in K_7$ the auxiliary word is $d(w)=(0-0)(1-0)(2-0)(3-2)(4-2)(5-3)(6-6)=0121220$ and the associated coefficient is 
	\[ c_w = \binom{0}{0}\binom{1}{0}\binom{2}{0}\binom{3}{2}\binom{2}{0}\binom{3}{1}\binom{3}{3} = 9.\]
	With these notations, we shall prove
	\begin{theo}\label{th:pbwnov} The map $\eta:U(A)\to A$ is explicitly given by
		\[ \eta(f_1\otimes\cdots \otimes f_i) = \sum_{w\in K_i} c_w \partial^{w_1}(f_1)\cdots\partial^{w_i}(f_i) \]
		for all  $f_1,\ldots,f_i\in A$, $i\ge1$.
	\end{theo}
	\begin{proof}
		We know from \eqref{eq:eta} that $\eta(f_1\otimes\cdots \otimes f_i) = \sum_{w\in K_i} d_w \partial^{w_1}(f_1)\cdots\partial^{w_i}(f_i)$, where $d_w$ is the number of maps $\phi$ in $Inc(i)$ such that $w_j=|\phi^{-1}(j)|$ for $2\leq j\leq i$.
		We shall prove $d_w=c_w$ by induction on the number of non zero coefficients $w_j$ as follows. 
		
		The claim is obvious when there is a single non zero coefficient (which is then necessarily $w_i=i-1$). Given $\phi\in Inc(i)$ with associated word $w=w_1\cdots w_i$, that is, $w_h:=|\phi^{-1}(h)|$ for all $1\le h\le i$, consider the smallest $j$ such that $w_j>0$ (necessarily, we also have $w_j\le j-1$). In particular, $j$ is the image under $\phi$ of $w_j$ integers that form an arbitrary subset $I_j$ of $[j-1]$, and moreover $\phi(k)>j$ for all $k\in[i-1]-I_j$.. 
		
		Using our previous notation, as our hypotheses imply $d_j+w_j=j-1$, there are 
		$\binom{d_j+w_j}{w_j}$ such subsets. Denote by $stand:[i-1]-I_j\to[i-1-w_j]$ the unique order preserving bijection, and define $\phi'\in Inc(i-w_j)$ according to $\phi'(k):=\phi\circ stand^{-1}(k)-w_j$. To check that $\phi'$ is indeed in $Inc(i-w_j)$, set $l:=max(I_j)$. If $k<l-w_j$, since $\phi\circ stand^{-1}(k)>j>l$ we have $\phi'(k)>l-w_j>k$, as expected. If instead $k\geq l-w_j$ then 
		$stand^{-1}(k)=k+w_j$, and therefore $\phi'(k)=\phi(k+w_j)-w_j>k+w_j-w_j=k$, as expected.
		
		Conversely, knowing the subset $I_j\subset[j-1]$ and the map $\phi'\in Inc(i-w_j)$ we can reconstruct $\phi$ according to $\phi(k)=j$ if $k\in I_j$ and $\phi(k)=\phi'\circ stand(k)+w_j$ if $k\in [i-1]-I_j$. 
		
		Finally, as $\phi$ had associated word $w=\overbrace{0\cdots0}^{j-1}w_j\cdots w_i$, it is easy to see that $\phi'$ has associated word $w'=\overbrace{0\cdots0}^{j-w_j}w_{j+1}\cdots w_i$, and conversely if we start with $I_j\subset[j-1]$ and $\phi'\in Inc(i-w_j)$ with associated word $w'$, then the resulting $\phi\in Inc(i)$, constructed as before from $I_j$ and $\phi'$, will have associated word $w$. Putting everything together, this shows that $d_w=\binom{d_j+w_j}{w_j}d_{w'}$, and the thesis follows from the inductive hypothesis applied to $d_{w'}$. 
	\end{proof}
	
	\begin{remark}\label{canorep} 
		In the construction we have just performed, when associating $c_w$ elements of $Inc(i)$ to $w\in K_i$, two choices of such elements are natural, namely taking for $I_j$ the subset $[w_j]$ or the subset $\{j-w_j,\dots,j-1\}$ of $[j-1]$. We call the second choice {\it canonical} for reasons that will become clear later in this section.
		In more abstract terms, to a map $\phi\in Inc(i)$ we can associate the sequence of values $\phi(1),\dots, \phi(i-1)$. The first choice we have just indicated corresponds to choosing the minimal element associated to $w$ in $Inc(i)$ for the lexicographical ordering, whereas the second choice amounts to choosing the maximal one.
	\end{remark}
	
	Theorem \ref{th:pbwnov} has, in view of Fla (\ref{pdfFlaPBW}) the following corollary, that we state as a theorem in view of its importance:
	\begin{theo}[PreLie PBW Theorem for Novikov algebras]\label{th:pbwinvnov} Let $N$ be a free Novikov algebra with generators $a^{[1]},\dots,a^{[i]}$. Then, we have
		in $\Q[N]$:
		
		$a^{[i]}\ast\dots \ast a^{[1]}=$\[\sum\limits_{k=1}^i\sum\limits_{I_1\coprod_<\dots\coprod_< I_k=[i]}\left(  \sum_{w\in K_{|I_1|}} c_w a^{[x_1^1]}_{w_1-1}\cdots a^{[x_{|I_1|}^1]}_{w_{|I_1|}-1} \right)\odot\dots\odot\left(  \sum_{z\in K_{|I_k|}} c_z a^{[x_1^k]}_{z_1-1}\cdots a^{[x_{|I_k|}^k]}_{z_{|I_k|}-1}  \right),\]
		with the notation $I_j=\{x_1^j,\dots,x_{|I_j|}^j\}$ and $I_1\coprod_<\dots\coprod_< I_k$ means that the blocks of the partition are ordered, in the sense that $min(I_1)<\dots <min(I_k)$.
	\end{theo}
	Applying this formula to a commutative algebra $A$ with a derivation and
	applying the evaluation map $ev$ to differential operators (Proposition \ref{mappingrule}) one gets the Leibniz type formula for $f_1,\dots,f_i$ in $A$:
	
	$f_1\partial\circ \dots \circ f_i\partial =$
	\[\sum\limits_{k=1}^i\sum\limits_{I_1\coprod_<\dots\coprod_< I_k=[i]}\left(  \sum_{w\in K_{|I_1|}} c_w \partial^{w_1}(f_{x_1^1})\cdots\partial^{w_{|I_1|}}(f_{x_{|I_1|}^1})  \right)\cdot\dots\cdot\left(  \sum_{z\in K_{|I_k|}} c_z \partial^{z_1}(f_{x_1^k})\cdots\partial^{z_{|I_k|}}(f_{x_{|I_k|}^k})  \right)\partial^i,\]
	
	Let us investigate now further the structure of the sets $K_i$.
	
	\begin{lemma}\label{freemagma} The binary operation
		\[-\vee-:K_i\times K_j\to K_{i+j},\qquad u=u_1\cdots u_i, v=v_1\cdots v_j\,\,\Rightarrow\,\, u\vee v := u_1\cdots u_iv_1\cdots v_{j-1}(v_j+1). \]
		makes $K=\coprod_{i\ge1} K_i$ a free magma, freely generated by $K_1$.
	\end{lemma}
	\begin{proof}
		
		First of all, it is easy to check that $u\in K_i,v\in K_j$ implies $u\vee v\in K_{i+j}$. In fact, if $1\le h \le i$ we have $u_1+\cdots+u_h\le h-1$ since $u\in K_i$, and if $1\le h<j$ we have $u_1+\cdots+u_i+v_1+\cdots+v_h=i-1+v_1+\cdots+v_h\le i+h-2$ since $v\in K_j$: finally, $u_1+\cdots+u_i+v_1+\cdots+v_j+1=(i-1)+(j-1)+1=i+j-1$, showing that  indeed $u\vee v\in K_{i+j}$.  
		
		Next, we need to check that given $i\ge2$ and $w=w_i\cdots w_i \in K_i$, there exists a unique factorization $w=u\vee v$. For this, denoting by $d(w)=d_1\cdots d_i$ the auxiliary word defined as before, we set
		\[j:=\max\{1\le h<i\,\,\mbox{s.t.}\,\, d_h = 0\},\]
		\[u= w_1\cdots w_j,\qquad v=w_{j+1}\cdots w_{i-1}(w_i-1).\]
		As $d_1=0$, $j$ is well-defined. Once we prove that $u\in K_j$, $v\in K_{i-j}$, it will be clear from the definitions that $w=u\vee v$, as desired. For $1\le h<j$ we have $w_1+\cdots+w_h\le h-1$ since $w\in K_i$: moreover, $d_j=0$ is equivalent to $w_1+\cdots+w_j=j-1$, thus $u=w_1\cdots w_j\in K_j$ as desired. For $1\le h <i-j$, by the way we defined $j$ we have that $d_{j+h}>0$ and thus $w_1+\cdots+w_j+w_{j+1}+\cdots+w_{j+h}=(j-1)+w_{j+1}+\cdots+w_{j+h}<j+h-1$, which in turn implies $w_{j+1}+\cdots+w_{j+h}\le h-1$: finally $i-1=w_1+\cdots+w_i=(j-1)+w_{j+1}+\cdots+w_{i}$ implies that $w_{j+1}+\cdots+w_{i}-1=i-j-1$, showing that $v=w_{j+1}\cdots w_{i-1}(w_i-1)\in K_{i-j}$, as desired.
		
		Finally, we need to show that if $w=u\vee v=u'\vee v'$ then $u=u'$ and $v=v'$. Suppose otherwise that $u\neq u'$: since both $u$ and $u'$ are left subfactors of $w$, either $u$ is a proper left subfactor of $u'$ or the other way around. Suppose to fix the ideas that $u=w_1\cdots w_j$ and $u'=w_1\cdots w_{j+h}$ for some $j,h\ge1$. Then $u\in K_j$, $u'\in K_{j+h}$ would imply that $j-1 = w_1+\cdots+w_j$ and $h+j-1=w_1+\cdots+w_{j+h}=(j-1)+w_{j+1}+\cdots+w_{j+h}$, which in turn would imply $w_{j+1}+\cdots+w_{j+h}=h$, which yields a contradiction since $w_{j+1}\cdots w_{j+h}$ is a left subfactor of $v=w_{j+1}\cdots w_{i-1}(w_i-1)$ and $v\in K_{i-j}$.
	\end{proof}
	
	As planar binary trees form a basis of the free magma on one generator, Lemma \ref{freemagma} implies the existence a bijection between $K_i$ and the set of planar binary rooted trees with $i$ leaves for all $i\ge1$. This can be described explicitly as follows. Given a planar binary rooted tree, we say that an edge is left pointing if the children is to the left of the parent, and right pointing otherwise. Notice that if the tree has $i$ leaves the number of left pointing edges is $i-1$ (and the same is true for right pointing edges). We define a correspondence $\phi$ from the left pointing edges to the leaves as follows: given a left pointing edge $e$, we begin by moving downwards  along $e$, and then we continue moving upwards, always picking the right pointing edge at each internal vertex we encounter, until we reach the leaf $\phi(e)$. For instance, in the following picture we colored every left pointing edge with the same color as the output leaf.
	\[T=
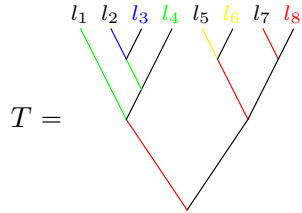
\xy {(-14,14)*{_{l_1}}};
	{(-10,14)*{_{l_2}}};
	{(-6,14)*{_{\color{blue}{l_3}}}};
	{(-2,14)*{_{\color{green}{l_4}}}};
	{(2,14)*{_{l_5}}};
	{(6,14)*{_{\color{yellow}{l_6}}}};
	{(10,14)*{_{l_7}}};
	{(14,14)*{_{\color{red}{l_8}}}};
	{\ar@{-}(0,-12);(8,0)};
	{\ar@{-}@[red](0,-12);(-8,0)};
	{\ar@{-}@[red](8,0);(4,8)};
	{\ar@{-}(8,0);(12,8)};
	{\ar@{-}@[green](-8,0);(-14,12)};
	{\ar@{-}(-8,0);(-6,4)};
	{\ar@{-}(-6,4);(-2,12)};
	{\ar@{-}@[green](-6,4);(-8,8)};
	{\ar@{-}@[blue](-8,8);(-10,12)};
	{\ar@{-}(-8,8);(-6,12)};
	{\ar@{-}@[yellow](4,8);(2,12)};
	{\ar@{-}(4,8);(6,12)};
	{\ar@{-}@[red](12,8);(10,12)};
	{\ar@{-}(12,8);(14,12)}; \endxy\]
	Finally, we order the leaves of $T$ from left to right $l_1,\ldots,l_i$, and we set $w_j:=|\phi^{-1}(l_j)|$ for all $1\le j\le i$, $w(T):=w_1\cdots w_i$. Then $T\to w(T)$ is the desired bijection between the set of planar binary rooted trees with $i$ leaves and $K_i$. For instance, for the tree in the above picture we have $w(T)=00120103$. 
	
	This construction can also be interpreted in terms of the correspondence between elements of $Inc(i)$ and $K_i$. 
	Let again $T$ be a planar binary rooted tree with leaves $l_1,\dots,l_i$, ordered from left to right. 
	Given a leaf at the extremity of a left pointing edge, proceed as above to map it to the leaf associated to the edge. Given instead a leaf at the extremity of a right pointing edge (other than the rightmost one), move downwards along right pointing edges till reaching a left pointing edge moving downwards. Move downwards along that edge and then upwards along right pointing edges till reaching a leaf. 
	This process defines a map from $\{l_1,\dots,l_{i-1}\}$ to $\{l_2,\dots,l_i\}$ and 
	a corresponding element of $Inc(i)$ that can be proven to be the canonical representative of the element of $K_i$ associated to the tree. In the example of the picture, the corresponding map in $Inc(8)$ is given by the sequence of values $(4,3,4,8,6,8,8)$.
	
	\section{The PBW theorem for Novikov algebras}\label{PBWc}
	
	We have used previously some general properties of the Poincar\'e-Birkhoff-Witt (PBW) theorem. The Theorem has been studied in detail in the context of Lie, pre-Lie algebras and their mutual relations in \cite{BS}, we investigate here the particular case of Novikov algebras. 
	
	Recall from \cite{Reutenauer} some generalities on the PBW theorem. Let $L$ be a Lie algebra and write $(U(L),\ast,\Delta)$ for an arbitrary enveloping algebra of $L$. The following map is then a coalgebra isomorphism:
	\[ \operatorname{pbw}: \Q[L]\to U(L) :x_1\odot\cdots\odot x_n\longmapsto \frac{1}{n!}\sum_{\sigma\in S_n} x_{\sigma(1)}\ast\cdots\ast x_{\sigma(n)}. \]
	This isomorphims induces a canonical decomposition of $U(L)=\bigoplus\limits_{n\in\N}U^{(n)}(L)$, the PBW decomposition, where, by definition
	$$U^{(n)}(L):=\operatorname{pbw}(\Q_n[L]),$$
	where we write $\Q_n[L]$ for the space of degree $n$ polynomials. The projection from $U(L)$ onto $U^{(n)}(L)$ orthogonally to the $U^{(j)}(L)$, $j\not= n$ is called the $n$-th canonical projection.
	
	When $N$ is a Novikov algebra, one can choose for $U(N)$, $\Q[N]$, the polynomial realization of the enveloping algebra, and the PBW map is then a linear endomorphism of $\Q[L]$ defined
	as 
	\begin{equation}\label{PBWr}
		x_1\odot\cdots\odot x_n\longmapsto \frac{1}{n!}\sum_{\sigma\in S_n} x_{\sigma(1)}\ast\cdots \ast x_{\sigma(n)}
	\end{equation}
	when $N$ is right Novikov and as
	\begin{equation}\label{PBWl}
		x_1\odot\cdots\odot x_n\longmapsto \frac{1}{n!}\sum_{\sigma\in S_n} x_{\sigma(1)}\ast^{op}\cdots \ast^{op} x_{\sigma(n)}
	\end{equation}
	when $N$ is left Novikov.
	\begin{lemma}
		The two PBW linear endomorphisms of $\Q[N]$ defined by Eq. (\ref{PBWr}) and Eq. (\ref{PBWl}) agree.
	\end{lemma}
	\begin{proof}
		Indeed, $x_{\sigma(1)}\ast\cdots \ast x_{\sigma(n)}=x_{\sigma(n)}\ast^{op}\cdots \ast^{op} x_{\sigma(1)}$, and the symmetric group is globally invariant under the map $(\sigma(1),\cdots ,\sigma(n))\longmapsto(\sigma(n),\cdots , \sigma(1))$.
	\end{proof}
	
	\begin{example}\label{expbwcalc2}
		For $n=2$, we get
		$$pbw(x\odot y)=\frac{1}{2}(x\ast y+y\ast x)=\frac{1}{2}(x\odot y+x\{y\}+y\odot x+y\{x\})$$
		$$=x\odot y+\frac{1}{2}(x\{y\}+y\{x\}),$$
		whereas 
		$$pbw(x\ast y)=pbw(x\odot y+x\{ y\})=x\odot y+\frac{1}{2}(3\cdot x\{y\}+y\{x\}).$$
		
		Accordingly, the PBW decomposition of $\Q[N]$ as an enveloping algebra for $N$ is non trivial:
		we get in the lowest degrees, $\Q^{(1)}[N]=N$ whereas, using polarization, $\Q^{(2)}[N]$ is the linear span of the elements $x\odot x+x\{x\}=x\ast x$ for $x\in N$.
		In particular the PBW decomposition of a degree 2 element in $\Q[N]$ reads:
		$$x\odot y = \left(- \frac{1}{2}(x\{y\}+y\{x\})\right) +\left(x\odot y+\frac{1}{2}(x\{y\}+y\{x\})\right).$$
		$$= \left(- \frac{1}{2}(x\{y\}+y\{x\})\right) +\frac{1}{2}(x\ast y+y\ast x).$$
	\end{example}
	
	\ \
	
	Recall now the definition of the Stirling numbers of the first kind $s(j,i)$: they are implicitly defined by
	$$x(x-1)\dots (x-j+1)=\sum_{i=1}^js(j,i)x^i.$$
	\begin{prop}\label{canproj}
		The $n$-th canonical projection is obtained as
		$$x_1\odot\cdots\odot x_m\longmapsto  \sum\limits_{k=n}^m\frac{s(k,n)}{k!}\sum\limits_{I_1\coprod\dots \coprod I_k=[m]}x^{\odot I_1}\ast\dots\ast x^{\odot I_k}$$
		or, equivalently,
		$$x_1\odot\cdots\odot x_m\longmapsto  \sum\limits_{k=n}^m\frac{s(k,n)}{k!}\sum\limits_{I_1\coprod\dots \coprod I_k=[m]}x^{\odot I_1}\ast^{op}\dots\ast^{op} x^{\odot I_k},$$
		where the blocks $I_i$ of the partitions of $[m]$ are non empty and $x^{\odot S}:=x_{s_1}\odot \dots\odot x_{s_p}$ if $S=\{s_1,\dots ,s_p\}$.
	\end{prop}
	
	\begin{proof}
		The Proposition is a direct application of a general Hopf algebraic formula for the canonical projections in an enveloping algebra, see \cite{patras1994algebre} or Thm 4.1.1 and 4.2.1 in \cite{Cartier2021}, taking into account the definition of the Hopf algebra structure on $(\Q[N],\ast,\Delta)$.
	\end{proof}
	\begin{example}
		Let us check the formula in the simplest non trivial case, checking that it agrees with Example \ref{expbwcalc2}:
		the first projection of $x\odot y$ is obtained as 
		$$s(1,1)x\odot y+\frac{s(2,1)}{2}(x\ast y+y\ast x)=x\odot y-\frac{1}{2}(x\ast y+y\ast x)=- \frac{1}{2}(x\{y\}+y\{x\}).$$
		The second projection is obtained as
		$$\frac{s(2,2)}{2}(x\ast y+y\ast x)=\frac{1}{2}(x\ast y+y\ast x).$$
	\end{example}
	
	When $N=(\mathcal A,\partial)$ is a commutative algebra with derivation and calculations in the enveloping algebra are interpreted as differential operators, we get that the first projection maps $x\odot y$ (corresponding to $xy\partial^2$) to the differential operator $- \frac{1}{2}(x\partial y+y\partial x)\cdot \partial$,
	whereas the second projection maps it to $\frac{1}{2}((x\partial)\circ( y\partial)+(y\partial)\circ( x\partial))$.
	More generally, we get:
	\begin{cor}
		In a commutative algebra with derivation, the following PBW decomposition holds:
		$$a_1\cdots a_m\partial^m=\sum\limits_{n=1}^m\left(\sum\limits_{k=n}^m\frac{s(k,n)}{k!}\sum\limits_{I_1\coprod\dots \coprod I_k=[m]}a_{I_1}\partial^{|I_1|}\circ\dots\circ a_{I_k}\partial^{|I_k|}\right).$$
	\end{cor}
	
	For reasons that should become clear in the next sections, where we will study constructions on pre-Lie algebras that have arisen from the theories of dynamical systems, control and numerical methods in analysis, this last decomposition should be meaningful in these application domains of differential calculus. However, we are not aware at the moment of its appearance elsewhere.

	\section{Pre-Lie exponential in Novikov algebras}\label{subsec:exp}
	
	Recall the definition of the pre-Lie exponential in an arbitrary pre-Lie algebra $L$:
	\begin{equation}\label{prelieexp}
		\exp_\triangleleft(x) = \sum\limits_{n=1}^\infty \frac{1}{n!}x^{\triangleleft n}= x+\frac{1}{2}x\triangleleft x +\frac{1}{6}(x\triangleleft x)\triangleleft x+\cdots 
	\end{equation}
	It is equivalently obtained as
	$$\exp_\triangleleft(x)=\log^\odot\circ\exp^\ast (x),$$
	where $\log^\odot$ is computed in $(\Q[L],\odot)$ and $\exp^\ast$ in $(\Q[L],\ast)$.

	We have chosen to treat here the case of right pre-Lie and Novikov algebras, the results hold {\it mutatis mutandis} for left pre-Lie and left Novikov algebras. We treat here the expression formally and do not investigate convergence issues.
	
	Pre-Lie exponentials are of the highest importance in view of applications of the theory of pre-Lie algebras. Besides being a key ingredient in the study of their structural properties \cite{BS,Kapr,Cartier2021} they appear for example in applications to dynamical systems and control theory \cite{AG} or in the relations between cumulants in noncommutative probability theories \cite{FoCM}.
	
	\ \
	
	Let us analyze first what can be learned on Novikov algebras from the theory of pre-Lie algebras..
	We denote by $\mathcal{T}_n$ the set of isomorphism classes of (non planar) rooted trees with $n$ vertices. We further denote by $\mathcal{T}_n^{inc}$ the set of isomorphism classes of pairs $(T,f)$, where $T\in\mathcal{T}_n$ and $f:V(T)\to\{1,\ldots,n\}$ is an order preserving bijection (we regard a rooted tree as the Hasse diagram of a poset structure on its set of vertices $V(T)$ as usual, putting the minimum at the root). It is known that the cardinality of $\mathcal{T}_n^{inc}$ is $(n-1)!$. We denote by $\Phi_n:\mathcal{T}_n^{inc}\to\mathcal{T}_n:(T,f)\to T$ the map forgetting the bijection $f$. Given a tree $T$ and a vertex $v\in V(T)$, let $f(v)$ be the number of children of $v$ (the \emph{fertility} of $v$). Notice that $\sum_{v\in V(T)} f(v)=n-1$: we denote by $\Psi_n:\mathcal{T}_{n}\to\mathcal{P}_{n-1}$ the map sending a tree $T$ to the corresponding partition of $n-1$. 
	
	Let $p\mathcal{L}(x)$ the free (right) pre-Lie algebra on one generator $x$: as we already recalled, it has a natural basis indexed by rooted trees, and a natural grading given by the number of vertices of the tree. We denote by $\bullet$ the generator of $p\mathcal{L}(x)$ (the tree with only one vertex), and by $\curvearrowleft$ the pre-Lie product to distinguish it from the pre-Lie product $\tl$ in the free Novikov algebra.
	
	The pre-Lie exponential $\exp_{\curvearrowleft}(\bullet)$ of the generator in the free pre-Lie algebra  $p\mathcal{L}(x)$ is well-known and particularly easy to compute. Indeed, as the product $T\curvearrowleft \bullet$ of a tree with the generator is obtained by all possible graftings of $\bullet$ on the vertices of $T$, it follows immediately (keeping track of the order of successive graftings) that the coefficient of a tree $T$ with $n$ vertices in the expansion of $x^{\triangleleft n}$ is equal to the number of increasing labelings of its vertices, that is:
	\[ \exp_{\curvearrowleft}(\bullet)=\sum_{n\ge1}\frac{1}{n!}\sum_{T\in\mathcal{T}^{inc}_n} \Phi_n(T)=\sum_{n\ge1}\frac{1}{n!}\sum_{T\in\mathcal{T}_n} N_TT,\]
	where $N_T=\#\Phi_n^{-1}(T)$ denotes the cardinality of the fiber of $\Phi_n:\mathcal{T}_n^{inc}\to\mathcal{T}_n$ at $T$.
	
	The pre-Lie exponential is usually re-expanded according to the equivalent formula
	\[ \exp_{\curvearrowleft}(\bullet) = \sum_{n\ge1}\sum_{T\in\mathcal{T}_n}\frac{1}{T!\sigma(T)} T, \]
	where $T!$ is the tree factorial\footnote{Given $v\in V(T)$ we denote by $d_v$ the number of descendants of $v$ (including $v$ itself). Then the tree factorial is defined as $T!:=\prod_{v\in V(T)}d_v$.}, and $\sigma(T)$ is the number of automorphisms of the rooted tree $T$. It is indeed well known that $\frac{n!}{T!}$ is the number of order preserving bijections $V(T)\to\{1,\ldots,n\}$
	so that $N_T:=\frac{n!}{T!\sigma(T)}$.
	
	As for Novikov algebras, the degree $n$ part of the free Novikov algebra $\mathcal{N}_n(x)$ has a natural basis indexed by $\mathcal{P}_{n-1}$. As the linear extension $\Psi:p\mathcal{L}(x)\to \mathcal{N}(x)$ of the maps $\Psi_n:\mathcal{T}_n\to\mathcal{P}_{n-1}$ is the unique morphism of pre-Lie algebras sending the generator $\bullet$ of $p\mathcal{L}(x)$ to the generator $x_{(0)}$ of $\mathcal{N}(x)$, we have
	\[\exp_\triangleleft(x_{(0)}) = \Psi\Big(\exp_{\curvearrowleft}(\bullet)\Big).\]
	
	In particular, this shows that
	\begin{prop}The pre-Lie exponential in Novikov algebras can be described as
		\[ \exp_{\triangleleft}(x_{(0)}) = \sum_{n\ge1}\frac{1}{n!}\sum_{p\in\mathcal{P}_{n-1}} \Big(\sum_{T\in\Psi_n^{-1}(p)}N_T\Big)x_p,   \]
		where the coefficient $\sum_{T\in\Psi_n^{-1}(p)}N_T$ in front of $x_p$ can be interpreted as the cardinality of the fiber of $\Psi_n\circ\Phi_n:\mathcal{T}_{n}^{inc}\to\mathcal{P}_{n-1}$ at $p$.
	\end{prop}
	
	To get a formula for the pre-Lie exponential in the free Novikov algebra $\mathcal N(x)$ better suited to Novikov algebras, we will follow three approaches. The third one will be discussed in relation to the calculation of the pre-Lie logarithm, and is postponed: see Proposition \ref{alterplexp}.
	
	The easiest solution is to use the same approach as for the study of the pre-Lie PBW theorem. 
	
	\begin{notation} We shall denote by $Leib(n)$ the set of maps from $[n]$ to $[n]$ such that $\phi(i)\leq i$ for all $i$ (these maps are called subexceedent functions in \cite{Mantaci01} but, according to the same article, this is not the standard meaning of the term). There is a natural map $\Theta_n:Leib(n)\to\mathcal{P}_n$ sending $\phi$ to the partition $n=\sum_{i=1}^n|\phi^{-1}(i)|$.
	\end{notation}
	
	To compute the degree $n$ component of the pre-Lie exponential, let us introduce elements $f_1,\dots,f_n$ in a free commutative algebra with a derivation $\partial$. For convenience we also introduce $n-1$ copies of $\partial$, denoted $\partial_1,\dots,\partial_{n-1}$. 
	By the Leibniz formula,
	$$(\dots ((f_1\tl f_2)\tl f_3 \dots )\tl f_n)=f_n\partial_{n-1}(f_{n-1}\partial_{n-2}(\dots f_2\partial_1(f_1))\dots )$$
	$$=\sum\limits_{\phi\in Leib(n-1)}f_nf_{n-1}^{(|\phi^{-1}(n-1)|)}\dots f_{1}^{(|\phi^{-1}(1)|)}.$$
	In words, a map $\phi$ corresponds to the term in the Leibniz formula such that $\partial_i$ acts as derivation on the function $f_{\phi(i)}$.
	
	This implies, using the dictionary between commutative algebras with a derivation and Novikov algebras and our notational conventions for elements in the free Novikov algebra, that:
	\begin{prop}\label{prop:expleib} In the free Novikov algebra $\mathcal N(x)$, we have:
		\begin{equation}\label{prelieexp2}
			\exp_\triangleleft(x_{-1})= \sum_{n\ge 1} \frac{1}{n!}\sum_{\phi\in Leib(n-1)}  x_{-1}x_{|\phi^{-1}(n-1)|-1}\dots x_{|\phi^{-1}(1))|-1},
		\end{equation}
		or, in the canonical basis of $\mathcal N(x)$,
		\begin{equation}\label{prelieexp3}
			\exp_\triangleleft(x_{(0)}) = \sum_{n\ge0} \frac{1}{(n+1)!}\sum_{p\in\mathcal{P}_n} N_p x_p,
		\end{equation}
		where $N_p$ is the cardinality of the fiber of $\Theta_{n}:Leib(n)\to\mathcal{P}_n$ at $p$.
	\end{prop}
	
	We shall discuss several combinatorial interpretations of the numbers $N_p$, starting from:
	
	\begin{lemma}\label{lem:formula}
		The coefficient $N_p$ of $x_p, p=(p_1,\dots,p_k)\in \mathcal{P}_{n,k}$ is equal to the number of (unordered) set partitions of $[n+k]$ with $k$ blocks of sizes $p_1+1,\dots,p_k+1$ such that the minimal elements of two distinct blocks are never two consecutive integers. 
	\end{lemma}
	
	\begin{proof}
		Let us show how to associate an element of $Leib(n)$ to such a set partition, the proof will follow (details are left to the reader). Consider a set partition $S_1\coprod\dots \coprod S_k$ of $[n+k]$ as in the Lemma, with therefore blocks of size at least 2, and ordered is such a way that $1=min(S_1)<\dots<min(S_k)$. We write $I$ for the set $\{min(S_1),\dots,min(S_k)\}$ and, as earlier, $stand$ for the standardization map (the unique increasing bijection) from $[n+k]-I$ to $[n]$. We have $i_1:=1=min(S_1)<i_2:=min(S_2)-1<\dots<i_k:=min(S_k)-(k-1)\leq n$.
		Let $\phi$ be defined, for $l\in[n]$, by:
		$$\phi(l):=i_j\iff stand^{-1}(l)\in S_j.$$
		This is well-defined, as $stand$ is surjective and the $S_j$ define a partition of $[n+k]$. To see that $\phi$ belongs to $Leib(n)$ we reason as follows. Given $p\in [n+k]-I$, there exists a unique $i$ such that  $min(S_i)<p<min(S_{i+1})$. In particular, $p\in S_j$ for a $j\leq i$, and $stand(p)=p-i$. We get: $\phi(p-i)=i_j=min(S_j)-j+1\leq min(S_i)-i+1\leq p-i$. This shows as desired that $\phi\in Leib(n)$, and it is also clear by the construction that $\Theta_n(\phi)=p$. Finally, it is easy to check that the construction is invertible at each step, so we can recover the partition $S^1\coprod\cdots\coprod S^k$ knowing the associated $\phi$.
	\end{proof}
	
	The previous combinatorial interpretation of $N_p$, together with an inclusion-exclusion argument, can be used to deduce a closed formula for $N_p$: to the best of our knowledge neither this nor other closed formulas for $N_p$ have appeared elsewhere in the literature. 
	
	First we set up some notations.
	\begin{notation} As customary, given $n_1+\cdots+n_k=n$, we shall denote by 
		\[{n\choose n_1,\ldots,n_k} =\frac{n!}{n_1!\cdots n_k!}\] 
    the associated multinomial coefficient. 
	\end{notation}
	\begin{notation} Given a partition $p=(p_1,\ldots,p_k)$ and a subset $I=\{i_1,\ldots,i_j\}\subset[k]$, we shall denote by $p_I:=p_{i_1}+\cdots+p_{i_j}$ and by
		\[{p_I\choose p_{i\in I}}:={p_I\choose p_{i_1},\ldots,p_{i_j}}.\] 
	\end{notation}
	\begin{notation} Given a partition $p=(p_1,\ldots,p_k)$, $p_1\ge\cdots\ge p_k\ge1$, if $k=a_1+\cdots+a_b$ with
		\[ p_1=\cdots =p_{a_1}>p_{a_1+1}=\cdots=p_{a_1+a_2}>\cdots>p_{a_1+\cdots+a_{b-1}+1}=\cdots=p_{k}\]
		we set $\sigma(p):=a_1!\cdots a_b!$ and call it the \emph{symmetry factor} of $p$.
	\end{notation}
	\begin{theo} For all partitions $p=(p_1,\ldots,p_k)$ we have
	\[ N_p=\frac{1}{\sigma(p)}\sum_{I_1\coprod_<\cdots\coprod_< I_h=[k]} (-1)^{k+h}|I_1|!\cdots|I_h|!{p_{I_1}\choose p_{i\in I_1}}\cdots{p_{I_h}\choose p_{i\in I_h}}{p_1+\cdots+p_k+h\choose p_{I_1}+1,\ldots,p_{I_h}+1}, \]
	where the sum runs over set partitions of $[k]$ ordered as in Theorem \ref{th:pbwinvnov}.
   \end{theo}
	\begin{proof} The desired formula will follow from Lemma \ref{lem:formula} and an inclusion-exclusion argument if we show that for a fixed set partition $I_1\amalg_<\cdots\amalg_< I_h=[k]$ the integer
		\[ |I_1|!\cdots|I_h|!{p_{I_1}\choose p_{i\in I_1}}\cdots{p_{I_h}\choose p_{i\in I_h}}{p_1+\cdots+p_k+h\choose p_{I_1}+1,\ldots,p_{I_h}+1} \]
		counts the number of ordered set partitions $J^1\amalg\cdots\amalg J^k=[p_1+\cdots+p_k+k]$ satisfying:
		\begin{enumerate} \item $|J^i|=p_i+1$ for all $i=1,\ldots,k$;
			\item for all $\ell=1,\ldots,h$ the set of minima $\{\min(J^i)\}_{i\in I_\ell}$ is a set of $|I_\ell|$ consecutive integers.
		\end{enumerate}
		This can be seen bijectively as follows. We start with an ordered set partition  $K^1\amalg\cdots\amalg K^h=[p_1+\cdots+p_k+h]$ into sets of cardinalities $|K^\ell|=p_{I_\ell}+1$, $\ell=1,\ldots,h$: the number of such partitions is counted by the multinomial coefficient ${p_1+\cdots+p_k+h\choose p_{I_1}+1,\ldots,p_{I_h}+1}$ in the formula. Let $\tau\in S_h$ be the permutation such that $1=m'_1:=\min(K^{\tau(1)})<\cdots<m'_h:=\min(K^{\tau(h)})$. Next, we define an order preserving injection $\phi:[p_1+\cdots+p_k+h]\hookrightarrow[p_1+\cdots+p_k+k]$ as follows: $\phi(1)=1$, and given $j>1$, if $\ell$ is the maximum among $1,\ldots,h$ such that $m'_\ell<j$, we set $\phi(j)=j+|I_{\tau(1)}|+\cdots+|I_{\tau(\ell)}|-\ell$. We denote by $K:=\operatorname{Im}(\phi)=\phi(K^1)\coprod\cdots\coprod\phi(K^h)\subset[p_1+\cdots+p_k+k]$ and by  $1=m_1:=\min(\phi(K^{\tau(1)}))<\cdots<m_h:=\min(\phi(K^{\tau(h)}))$: notice that, by construction, for all $\ell=1,\ldots,h$ the difference between $m_\ell=\phi(m'_\ell)$ and $\phi(m'_\ell+1)$ (the next smaller entry in $K$) is precisely $|I_{\tau(\ell)}|-1$. In particular, 
		\begin{equation*} ([p_1+\cdots+p_k+k]-K)\coprod\{m_1,\ldots,m_h\} = \coprod_{\ell=1}^h\{m_\ell,m_\ell+1,\ldots,m_\ell+|I_{\tau(\ell)}|-1\}.
		\end{equation*} Now, for all $\ell=1,\ldots,h$, if $I_{\tau(\ell)}=\{a_1,\ldots,a_b\}$, we choose an ordered partition $\phi(K^{\tau(\ell)})-\{m_\ell\} = \overline{J}^{a_1}\coprod\cdots\coprod \overline{J}^{a_b}$ into sets of cardinalities $|\overline{J}^{a_r}|=p_{a_r}$, $r=1,\ldots,b$, together with a permutation $\rho\in S_{b}$: counting the number of possible choices gives the factor $|I_{\tau(\ell)}|!{p_{I_{\tau(\ell)}}\choose p_{i\in I_{\tau(\ell)}} }$ in the formula. Finally, we set $J^{a_r}:=\{m_\ell+\rho(r)-1\}\coprod\overline{J}^{a_r}$ for all $r=1,\ldots,b$. This gives the desired partition $J^1\coprod\cdots\coprod J^k=[{p_1+\cdots+p_k+k}]$ satisfying conditions 1 and 2 above.
	\end{proof}
	\begin{example} For instance, if $p=(p_1,p_2)$ we get
	\[ \sigma(p)N_p={p_1+p_2+2\choose p_1+1,p_2+1}-2{p_1+p_2\choose p_1,p_2},\]
	and if $p=(p_1,p_2,p_3)$ we get
	\begin{eqnarray} \nonumber \sigma(p)N_p&=&{p_1+p_2+p_3+3\choose p_1+1,p_2+1,p_3+1}-2{p_1+p_2\choose p_1,p_2}{p_1+p_2+p_3+2\choose p_1+p_2+1,p_3+1}\\\nonumber&&-2{p_1+p_3\choose p_1,p_3}{p_1+p_2+p_3+2\choose p_1+p_3+1,p_2+1}-2{p_2+p_3\choose p_2,p_3}{p_1+p_2+p_3+2\choose p_2+p_3+1,p_1+1}\\\nonumber&&+6{p_1+p_2+p_3\choose p_1,p_2,p_3}
	\end{eqnarray}
\end{example}
	
	Our next aim is to show that the numbers $N_p$ can also be understood combinatorially as a statistics on permutations. The fact that such a connection between a pre-Lie exponential formula and symmetric groups exists is maybe not so surprising: it is well-known that the classical Baker-Campbell-Hausdorff (BCH) formula is closely related to the inverse of the pre-Lie exponential (see \cite{EFP14}), and that the BCH formula can be expressed in terms of descent classes in symmetric groups \cite[Chap. 3]{Reutenauer}, that is the families of permutations $\sigma\in S_n$ with the same descent set 
	\[desc(\sigma)=\{1\leq i<n \ s.t.\  \sigma(i)>\sigma(i+1)\}.\]
	
	\begin{notation} We shall denote by $Leh(n)$ the set of words $l_1\cdots l_n\in\mathbb{N}^{\times n}$ such that $l_i\le n-i$ for all $1\le i\le n$. We call a word in $Leh(n)$ a \emph{Lehmer code}.
	\end{notation}
	
	Given a permutation $\sigma\in S_n$, its associated \emph{Lehmer code} is the word $\ell(\sigma)=l_1\cdots l_n\in\mathbb{N}^{\times n}$ defined by
	\[ l_i:=\#\{j>i\,\,\mbox{s.t.}\,\, \sigma(j)<\sigma(i) \} \]
	Representing $\sigma$ with the word $\sigma(1)\cdots\sigma(n)$, then $l_i$ is the number of letters to the right of $\sigma(i)$ which are smaller than $\sigma(i)$. It is a well known fact that this construction establishes a bijective correspondence $\ell:S_n\to Leh(n)$.	\begin{lemma}\label{lehcodbij}
		The map $Code:\ Leib(n)\longmapsto Leh(n)$  defined by
		$$Code(\phi):=(\phi(n)-1,\phi(n-1)-1,\dots, \phi(1)-1)$$
		is a bijection. 
	\end{lemma}
	
	The Lemma immediately follows from the definitions. It suggests to associate a partition $L(\sigma)\in\mathcal{P}_n$ to $\sigma\in S_n$ by considering the multiplicities in the Lehmer code $\ell(\sigma)$. For instance, for the following permutations in $S_9$
	\[ \sigma\:=\qquad 362857194,\qquad 812374596,\qquad 648327915 \]
	the Lehmer codes are
	\[ \ell(\sigma)\:=\qquad 241422010,\qquad 700030010,\qquad 535212200\]
	and the corresponding partitions are
	\[ L(\sigma):=\qquad (3,2,2,2),\qquad (6,1,1,1),\qquad(3,2,2,1,1). \]
	
	Finally, given $p\in\mathcal{P}_n$ the Lemma \ref{lehcodbij} implies that
	\[ N_p= \#L^{-1}(p) \]
	the cardinality of the fiber of $L$ over $p$.
	Together with Proposition \ref{prop:expleib}, this shows:
	
	\begin{theo}\label{preLexp} The pre-Lie exponential 
		of the generator $x_{(0)}$ in the free Novikov algebra $\mathcal{N}(x)$ is governed by the statistics defined by the map $L$ on permutations. That is, it is given by the formula
		\[ \exp_\triangleleft(x_{(0)})= \sum_{n\ge0} \frac{1}{(n+1)!}\sum_{\sigma\in S_n} x_{L(\sigma)}. \]
	\end{theo}
	
	\begin{remark} The links with descent classes of permutations to which we alluded can be made more precise. It is a direct consequence of the definition of the Lehmer code of a permutation $\sigma\in S_n$ that $\sigma(i)>\sigma(i+1)\iff l_{i+1}< l_i$ or, equivalently,  $\sigma(i)<\sigma(i+1)\iff l_{i+1}\geq l_i$.
		Therefore, defining the set of descents of a Lehmer code $l_1...l_n$ as $\{i\leq n-1|l_{i+1}<l_i\}$, the bijection between Lehmer codes and permutations restricts to bijections between Lehmer codes with a fixed set of descents $S$ and permutations with the same descent set $S$.
	\end{remark}
	
	\begin{remark}\label{alternat} Another link exists between Lehmer codes and descents. 
		Denoting by $\mathcal{P}_{n,k}$ the set of partitions of $n$ of length $k$, and by $a_{n,k}$ the Eulerian numbers, that is, the number of permutations $\sigma\in S_n$ with $d(\sigma)=k$, where $d(\sigma)$ is the descent number \[d(\sigma):=\#\{1\le i<n \mbox{ s.t. } \sigma(i)>\sigma(i+1)\},\] we have
		\[ \sum_{p\in\mathcal{P}_{n,k}} N_p = a_{n,k-1}. \]
		The result is attributed in \cite{Mantaci01} to D. Dumont (unpublished). The same article contains an independent bijective proof of the formula and introduce another bijection between elements of $Leib(n)$ and permutations than the one we use. The authors, Mantaci and Rakotondrajao, also point out that ``their coding seems to be appropriate to
		transfer to the subexceedant functions (the term they use for elements of the $Leib(n)$) some properties that are typical of permutations, namely the distribution of
		certain statistics called Eulerian". It is not clear to us at the moment whether this alternative coding might be useful to further investigate the pre-Lie structure of Novikov algebras.
	\end{remark}
	
	Finally, we directly relate the Lehmer code approach to the pre-Lie exponential for Novikov algebras to the classical pre-Lie approach by
	explicitly constructing a
	bijection $\Gamma_n:S_{n-1}\to\mathcal{T}_{n}^{inc}$ making the following diagram commutative
	\[ \xymatrix{ \mathcal{T}_n^{inc}\ar[r]^-{\Phi_n} & \mathcal{T}_n\ar[d]^-{\Psi_n} \\ S_{n-1}\ar[u]^-{\Gamma_n}\ar[r]_-L & \mathcal{P}_{n-1} }
	\]
	where $L:S_{n-1}\to\mathcal{P}_{n-1}$ is the map introduced in the previous subsection counting the multiplicities in the Lehmer code. 
	
	The existence of bijections between permutations and increasing rooted trees is classical since, at least, the work of X. Viennot. We define the one, denoted $\Gamma_n$, with the desired property recursively as follows. For $n=1$ all the involved sets are by definition singletons. Assume inductively we have defined $\Gamma_n$. To define $\Gamma_{n+1}$, we represent a permutation $\sigma\in S_{n}$ as the word $\sigma(1)\cdots\sigma(n)$, and let $\sigma'\in S_{n-1}$ be the permutation obtained by removing the letter $n$ from $\sigma$. Number the positions in which you could insert a new letter $n$ in $\sigma'$ from $1$ to $n$, with $1$ corresponding to inserting $n$ in the rightmost position and $n$ corresponding to inserting $n$ in the leftmost position. Write $\sigma = \sigma'\vee_i n$ if $\sigma$ is obtained from $\sigma'$ by inserting $n$ in the $i$-th position. Also, given a tree $T'\in\mathcal{T}_n^{inc}$, denote by $T'\swarrow_i \bullet_{n+1}\in\mathcal{T}_{n+1}^{inc}$ the tree obtained from $T'$ by attaching a new leaf labeled by $n+1$ to the vertex labeled by $i$ in $T'$. Finally, we define $\Gamma_{n+1}:S_n\to\mathcal{T}_{n+1}^{inc}$ according to
	\[ \sigma = \sigma'\vee_i n\qquad\Rightarrow\qquad\Gamma_{n+1}(\sigma) := \Gamma_n(\sigma')\swarrow_i\bullet_{n+1}\]
	
	The recursion starts sending the unique element in $S_0:=\{0\}$ to the tree $\Gamma_1(0)=\bullet_1$. For the unique element $1\in S_1$, we have $\Gamma_2(1) = \bullet_1\swarrow_1\bullet_2=\xy {(0,-4.5)*{_1}};{(0,4.5)*{_2}};{\ar@{-}(0,-2)*{\bullet};(0,2)*{\bullet}}\endxy$. As a more complicated example, for the permutation
	\[ 35412 = 3412\vee_4 5 = (312\vee_3 4)\vee_4 5 = ((12\vee_3 3)\vee_3 4)\vee_4 5= (((1\vee_1 2)\vee_3 3)\vee_3 4)\vee_4 5  \]
	the corresponding tree is
	\[ \Gamma_6(35412) = ((((\bullet_1\swarrow_1\bullet_2)\swarrow_1\bullet_3)\swarrow_3\bullet_4)\swarrow_3\bullet_5)\swarrow_4\bullet_6= \xy 
	{(0,-6)*{_1}};{(-6,0)*{_2}};{(6,0)*{_3}};
	{(-2,4)*{_4}};{(10,4)*{_5}};{(0,11.2)*{_6}};
	{\ar@{-}(0,-4)*{\bullet};(-4,0)*{\bullet}};
	{\ar@{-}(0,-4)*{\bullet};(4,0)*{\bullet}};
	{\ar@{-}(4,0)*{\bullet};(0,4)*{\bullet}};
	{\ar@{-}(4,0)*{\bullet};(8,4)*{\bullet}};
	{\ar@{-}(0,4)*{\bullet};(0,9)*{\bullet}};
	
	\endxy\]
	It is not hard to check that $\Gamma_{n+1}$ is a bijection making the required diagram commutative. In fact one can show the following more precise statement: the number of children of the vertex labeled by $i$ in the tree $T(\sigma)$ coincides with the multiplicity of the letter $i-1$ in the Lehmer code of $\sigma$. 
	For this, notice that 
	\[ \sigma = \sigma'\vee_i n\qquad\Rightarrow \qquad \ell(\sigma) = \ell(\sigma')\vee_i (i-1) \]that is, if $\sigma$ is obtained from $\sigma'$ by inserting the letter $n$ in the $i$-th position, then $\ell(\sigma)$ is obtained from $\ell(\sigma')$ by inserting the letter $i-1$ in the same position: together with an easy induction, this implies the claim.

	\section{The pre-Lie logarithm} \label{sec:lie} 
	Recall the definition of the pre-Lie logarithm, also known as Magnus map or, in the context of differential equations, Magnus expansion (as it provides a way to implicitly or numerically determine the logarithm of the solution of a matrix or operator linear differential equation), see e.g. \cite[Chap. 6]{Cartier2021} for a survey. Given a Novikov algebra $N$, the pre-Lie logarithm is the set automorphism of $N$ given for $x\in N$ by:
	\begin{equation}\label{prelielog} 
		\log_\triangleleft(x) = \log_\ast\circ\exp_\odot (x),
	\end{equation}
	where $\log_\ast$ is the logarithm computed in $(\Q[N],\ast)$ and $\exp_\odot$ the exponential computed in $(\Q[N],\odot)$.
	We get (formally, as we do not discuss convergence issues here):
	\begin{theo}\label{theoflaplog}
		For $f\in N$, we have:
		\begin{equation}
			\log_\triangleleft(f) = 
			f+\sum\limits_{n=2}^\infty \sum\limits_{i=1}^{n-1}\frac{(-1)^i}{i+1}\left(\sum\limits_{k_1+\dots+k_i=n-1}\frac{1}{k_1!\dots k_i!}(\dots (f\{f^{\odot k_1}\})\dots )\{f^{\odot k_i}\}\right),
		\end{equation}
		where the $k_i$ are positive integers. 
		When $N=(\mathcal A,\partial)$ is furthermore a commutative algebra with derivation we obtain:
		\begin{equation}
			\log_\triangleleft(f) = f+\sum\limits_{n=2}^\infty \sum\limits_{i=1}^{n-1}\frac{(-1)^i}{i+1}\left(\sum\limits_{k_1+\dots+k_i=n-1}\frac{1}{k_1!\dots k_i!}f^{k_i}\partial^{k_i}(\dots(f^{k_1}\partial^{k_1}(f))\dots)\right).
		\end{equation} 
	\end{theo}
	\begin{proof}We assume here that the reader is familiar with standard arguments in Hopf algebras (see \cite{Cartier2021}, in particular Remark 5.2.2 on the Baker-Campbell-Hausdorff formula, a topic closely connected to the ones we are addressing in this article).
		We have,
		$$\log_\triangleleft(f) = \log_\ast(\exp_\odot (f)).$$
		However, as $\exp_\odot (f)$ is the exponential of a primitive element in the enveloping algebra, it is a group-like element. This implies that the action of $\log_\ast$ on  $\exp_\odot (f)$ identifies with the action of the first canonical projection. 
		Using Proposition \ref{canproj}, we get
		$$\log_\triangleleft(f) =\sum\limits_{n=1}^\infty\frac{1}{n!} \left(\sum\limits_{j=1}^{n}\frac{(-1)^{j-1}}{j}\left(\sum\limits_{I_1\coprod\dots\coprod I_j=[n]}f^{\odot |I_1|}\ast\dots\ast f^{\odot |I_j|}\right)\right).$$
		As $\log_\triangleleft(f)$ belongs to $N$, non linear components have to vanish, and the expression simplifies to
		$$\log_\triangleleft(f)=\sum\limits_{n=1}^\infty \frac{1}{n!}\left(\sum\limits_{j=1}^{n}\frac{(-1)^{j-1}}{j}\left(\sum\limits_{I_1\coprod\dots\coprod I_j=[n]\atop |I_1|=1}f\{f^{\odot |I_2|}\ast \dots\ast f^{\odot |I_j|}\}\right)\right)$$
		$$=\sum\limits_{n=1}^\infty \frac{1}{n!}\left(\sum\limits_{j=1}^{n}\frac{(-1)^{j-1}}{j}\left(\sum\limits_{1+p_2+\dots +p_j=n}{n\choose 1\ p_2\ \dots \ p_j}f\{f^{\odot p_2}\ast \dots\ast f^{\odot p_j}\}\right)\right).$$
		The first part of the Theorem follows using Equations (\ref{sb1}) and (\ref{envprod}).
		The last statement follows using the mapping rule of Proposition \ref{mappingrule}.
	\end{proof}
	
	Let us consider first the lessons that can be learned from pre-Lie algebras. In the free pre-Lie algebra, the pre-Lie logarithm is obtained as  (\cite{Mur,CEM,BS,CP}):
	\begin{prop}
		\begin{equation}
			\log_\curvearrowleft(\bullet) = \sum\limits_T\frac{\omega(T)}{\sigma(T)}T,
		\end{equation}
		where $T$ runs over (non planar) trees,
		\begin{equation}
			\omega(T)=\sum\limits_{k=1}^{|T|}\frac{(-1)^{k-1}}{k}\omega_k(T),
		\end{equation}
		and $\omega_k(T)$ stands for the number of surjective strictly increasing maps from $T$ to $[k]$.
	\end{prop}
	We get, in the free Novikov algebra:
	\begin{cor}
		\begin{equation}
			\log_\triangleleft(x_{(0)}) = \sum\limits_p\left(\sum\limits_{\Psi(T)=p}\frac{\omega(T)}{\sigma(T)}\right)x_p,
		\end{equation}
		where $p$ runs over integer partitions and we recall that $\Psi$ stands from the canonical map from trees to partitions counting the arities of vertices with multiplicities.
	\end{cor}
	
	We will aim now in this section at finding a representation of $\log_\triangleleft(x_{(0)})$ better suited to Novikov algebras, looking for an expression of the coefficients 
	$n_{p,k}$:
	\begin{equation}\label{flaNovLlog}
		\log_\triangleleft(x) = x+\sum\limits_{m=2}^\infty\sum\limits_{p\in\mathcal{P}_{m-1}}\left(\sum\limits_{i=2}^{m}\frac{(-1)^{i-1}}{i}n_{p,i}\right)\frac{x_p}{p!},
	\end{equation}
	where, for notational tractability we abbreviate $x_{(0)}$ to $x$. We also set $n_{(0),1}:=1$ to account for the term $x$ in the expansion.
	By Theorem
	\ref{theoflaplog}, calculating $n_{p,i+1}$ boils down to computing the coefficient of $\frac{x_p}{p!}$ in the expansion of 
\begin{equation}\label{eq:npi}
	\sum\limits_{k_1+\dots+k_i=m-1}\frac{1}{k_1!\dots k_i!}(\dots (x\{x^{\odot k_1}\})\dots )\{x^{\odot k_i}\}.
\end{equation}	
	We will do so by a recursion on $i$.
	\begin{notation}
		Let $p',p$ be two partitions such that $|p|=|p'|+k$. We write $C_{p,p'}$ for the coefficient of $\frac{k!}{p!}x_p$ in the expansion of
		$\frac{1}{p'!}x_{p'}\{x^{\odot k}\}$. 
	\end{notation}
	
	Using the algebra with derivation realization of the free Novikov algebra, $C_{p,p'}$ also computes the coefficient of $\frac{k!}{p!}x_p$ in the expansion of $\frac{1}{p'!}x_1^k\partial^k(x_{p'_1-1}\dots x_{p'_{l(p')}-1}x_{-1}^{|p'|-l(p')+1})$. This implies in particular that $C_{p,p'}=0$ unless $p'<p$ and $l(p)\le|p'|+1$.

	\begin{notation}
		Let now $p',p$ be two partitions such that $p'<p$. We say that the pair is {\em admissible} if and only if $l(p)\leq |p'|+1$ and write then $p'<_{a}p$. \emph{Caution:} notice that $<_{a}$ is not an order relation (it is not transitive).
	\end{notation}
	
	With these notations set up, we now see that 
	\begin{eqnarray}\nonumber \sum_{p\in\mathcal{P}_{m-1}} n_{p,i+1}\frac{x_p}{p!} &=& \sum\limits_{k_1+\dots+k_i=m-1}\frac{1}{k_1!\dots k_i!}(\dots (x\{x^{\odot k_1}\})\dots )\{x^{\odot k_i}\}\\\nonumber &=&\sum_{k_i=1}^{m-1}\frac{1}{k_i!}\left(\sum\limits_{k_1+\dots+k_{i-1}=m-k_i-1}\frac{1}{k_1!\dots k_{i-1}!}(\dots (x\{x^{\odot k_1}\})\dots )\{x^{\odot k_{i-1}}\} \right)\{x^{\odot k_i}\} \\\nonumber&=& \sum_{k=1}^{m-1}\frac{1}{k!}\left(\sum_{p'\in\mathcal{P}_{m-k-1}}n_{p',i}\frac{x_{p'}}{p'!}\right)\{x^{\odot k}\}\\\nonumber &=&\sum_{k=1}^{m-1}\sum_{p'\in\mathcal{P}_{m-k-1}}\frac{1}{k!}n_{p',i}\left(\sum_{p\in\mathcal{P}_{m-1}} C_{p,p'}k!\frac{x_p}{p!}\right)\\\nonumber&=&\sum_{p\in\mathcal{P}_{m-1}}\left(\sum_{p'<_a p} C_{p,p'}n_{p',i}\right)\frac{x_p}{p!}
	\end{eqnarray}
	This produces the desired recursion 
		$$n_{p,i+1}=\sum\limits_{p'<_ap}C_{p,p'}n_{p',i}$$
	for computing the coefficients $n_{p,k}$. The basis of the recursion is $n_{p,1}=\delta^p_{(0)}$ (Kroenecker's delta). Our next aim is to have a better understanding of the coefficients $C_{p,p'}$ appearing in the recursion. 
	\begin{lemma} For $(p,p')$ admissible,
		consider the set $M(p,p')$ of $2\times(|p'|+1)$ matrices $A$ such that:
		\begin{itemize} \item the non-zero entries in the first row coincide with the non-zero entries of $p$, and the non-zero entries of the second row coincide with the non-zero entries of $p'$.
			
			\item $a_{2,i}\ge a_{2,i+1}$ for all $1\le i\le |p'|$.
			
			\item $a_{1,i}\ge a_{2,i}$ for all $1\le i\le|p'|+1$.
		\end{itemize}
		Then we have
		\[ C_{p,p'} = \sum_{A\in M(p,p')}\prod_{i=1}^{|p'|+1}\binom{a_{1,i}}{a_{2,i}}.\]
	\end{lemma}
	
	\begin{proof} Set $l:=l(p')$ and consider the coefficient of
		$x_p$ in the expansion of $x_{-1}^k\partial^k(x_{p'_1-1}\dots x_{p'_{l}-1}x_{-1}^{|p'|-l+1})$. Recall that $\partial x_i=x_{i+1}$. By the Taylor formula, $\partial^k(x_{p'_1-1}\dots x_{p'_{l}-1}x_{-1}^{|p'|-l+1})$ is the sum
		$$\sum\limits_{i_1+\dots +i_{|p'|+1}=k}{k\choose{i_1,\dots ,i_{|p'|+1}}} x_{p'_1+i_1-1}\dots x_{p'_{l}+i_{l}-1}x_{i_{l+1}-1}\dots x_{i_{|p'|+1}-1}.$$
		Each term of this sum can be associated to a matrix 
		
		\[ \left(\begin{array}{cccccc} 
			p'_1+i_1&\dots &p'_{l}+i_{l}&i_{l+1}&\dots &i_{|p'|+1} \\ 
			p'_1&\dots &p'_{l}&0&\dots &0
		\end{array}\right)\]
		and the terms that contribute to $x_p$ are exactly those terms associated to matrices in $M(p,p')$.
		
		The corresponding combinatorial coefficient contributing to the calculation of $C(p,p')$ is:
		$$\frac{p!}{k!p'!}{k\choose{i_1,\dots ,i_{|p'|+1}}}=\frac{(p'_1+i_1)!\dots (p'_{l}+i_{l})!i_{l+1}!\dots i_{|p'|+1}!}{(p'_1!i_1!)\dots (p'_l!i_l!)i_{l+1}!\dots i_{|p'|+1}!}$$
		$$={p'_1+i_1\choose p'_1}\dots {p'_{l}+i_{l}\choose p'_l}{i_{l+1}\choose 0}\dots {i_{|p'|+1}\choose 0},$$
		which concludes the proof.
	\end{proof}

	For instance, if $p=(2,2,1)$ and $p'=(2,1)$ the set $M(p,p')$ consists of the following matrices
	\[ \left(\begin{array}{cccc} 2 & 2 & 1 & 0 \\ 2 & 1 & 0 & 0
	\end{array}\right)\qquad\left(\begin{array}{cccc} 2 & 2 & 0 & 1 \\ 2 & 1 & 0 & 0
	\end{array}\right)\qquad\left(\begin{array}{cccc} 2 & 1 & 2 & 0 \\ 2 & 1 & 0 & 0
	\end{array}\right)\qquad\left(\begin{array}{cccc} 2 & 1 & 0 & 2 \\ 2 & 1 & 0 & 0
	\end{array}\right)\]
	and thus
	\[C_{(2,2,1),(2,1)}= \binom{2}{2}\binom{2}{1}\binom{1}{0}\binom{0}{0}+\binom{2}{2}\binom{2}{1}\binom{0}{0}\binom{1}{0}+\binom{2}{2}\binom{1}{1}\binom{2}{0}\binom{0}{0}+\binom{2}{2}\binom{1}{1}\binom{0}{0}\binom{2}{0} = 6. \]
	
	Finally, putting the previous results together, we get the Theorem:
	\begin{theo}
		We have for the coefficients $n_{p,k}$ of the pre-Lie logarithm in the free Novikov algebra
		$$\log_\triangleleft(x) = x+\sum\limits_{n=2}^\infty\sum\limits_{p\in\mathcal{P}_{n-1}}\left(\sum\limits_{i=2}^{|p|+1}\frac{(-1)^{i-1}}{i}n_{p,i}\right)\frac{x_p}{p!}$$
		that
		$$n_{p,k+1}=\sum\limits_{(0)=p_0<_a \dots <_a p_k=p}\prod\limits_{h=1}^kC_{p_{h},p_{h-1}}.$$
	\end{theo}
	\begin{remark}\label{rem:combnpi} We can interpret the integers $n_{p,i}\in\mathbb{N}$ combinatorially as counting the numbers of pairs $(T,\ell)$, where $T$ is a planar rooted tree with associated partition $\Psi(T)=p$, and where $\ell:V(T)\to[1,i]$ is a surjective and strictly monotone map. This follows from the forthcoming Lemma \ref{lem:combQ} and a standard inclusion-exclusion argument.
	\end{remark}
	As partitions can be represented by tableaux and pairs $(p',p)$ with $p'<p$ by skew tableaux,
	we can organize the computation of the coefficients $n_{p,i}$ in terms of tableaux of shape $p$ satisfying certain rules. Namely, given an ordered sequence $(0)=p_0< \dots <_a p_k=p$ as above, it is encoded by the Young tableau associated to $p$ with a sequence of sub tableaux associated to $p_1,\dots, p_k$. We represent this pictorially by putting the integer $h$ in all boxes of the skew tableau $p_h-p_{h-1}$, for all $h=1,\ldots,k$: then the fact that all successive pairs have to be admissible imposes constraints on the resulting decoration of $p$.
	
	More precisely, denote by $T_{p,k}$ the set of maps $s$ from the boxes of $p$ to $[k]$ satisfying:
	\begin{itemize} \item $s$ is surjective; \item $s$ is non-decreasing along the rows and columns of $p$; \item for all $1\le h\le k$, if $N_{<h}=\#s^{-1}(\{1,\dots,h-1\})$ is the number of boxes with labels $<h$ (in particular $N_{<1}=0$), then the boxes labeled by $h$ are contained in the first $N_{<h}+1$ rows of $p$ (for instance, the boxes labeled by $1$ are all contained in the first row of $p$, and if there are $i$ such boxes then the ones labeled by $2$ are all contained in the first $i+1$ rows of $p$, and so on).
	\end{itemize} For each $s\in T_{p,k}$ and $0\le h \le k$, denote by $s_{\le h}$ the tableau obtained from $s$ by removing the boxes labeled $h+1,\ldots,k$ (in particular $s_{\le0}$ is the empty tableau), as well as the underlying partition. 
	Then we have the following combinatorial formula for $n_{p,k+1}$: as already remarked $n_{p,1}=\delta^p_{(0)}$, while for $k\ge1$
	\begin{equation}\label{tablpllog} n_{p,k+1} = \sum_{s\in T_{p,k}} \prod_{h=1}^{k} C_{s_{\le h},s_{\le h-1}}.\end{equation}
	
	For instance, when $p=(2,1,1)$ the set $T_{p,4}$ consists of the standard tableaux
	\[\begin{array}{cc} 1 & 2 \\ 3 &\\ 4&
	\end{array}\qquad\begin{array}{cc} 1 & 3 \\ 2 & \\4&
	\end{array}\qquad\begin{array}{cc} 1 & 4 \\ 2 & \\3&
	\end{array} \]
	from which we get
\begin{eqnarray}\nonumber n_{(2,1,1),5} &=& C_{(2,1,1),(2,1)}C_{(2,1),(2)}C_{(2),(1)}+C_{(2,1,1),(2,1)}C_{(2,1),(1,1)}C_{(1,1),(1)}+C_{(2,1,1),(1,1,1)}C_{(1,1,1),(1,1)}C_{(1,1),(1)} \\\nonumber &=& 2 \cdot 2\cdot 2+ 2 \cdot 4\cdot 1+6\cdot1\cdot1 = 22 
\end{eqnarray}

	The set $T_{p,3}$ consists of the tableaux
	\[\begin{array}{cc} 1 & 1 \\ 2 &\\ 3&
	\end{array}\qquad\begin{array}{cc} 1 & 2 \\ 2 &\\ 3&
	\end{array}\qquad \begin{array}{cc} 1 & 2 \\ 3 &\\ 3&
	\end{array}\qquad\begin{array}{cc} 1 & 3 \\ 2 &\\ 3&
	\end{array} \]
	from which we get
	\begin{eqnarray}\nonumber n_{(2,1,1),4} &=& C_{(2,1,1),(2,1)}C_{(2,1),(2)}+C_{(2,1,1),(2,1)}C_{(2,1),(1)}+C_{(2,1,1),(2)}C_{(2),(1)}+C_{(2,1,1),(1,1)}C_{(1,1),(1)}\\\nonumber&=& 2\cdot 2+ 2\cdot 3 + 1\cdot 2 + 5\cdot 1 = 17
	\end{eqnarray}
	The set $T_{p,2}$ consists of the single tableau
	\[\begin{array}{cc} 1 & 1 \\ 2 & \\2&
	\end{array} \]
	from which we get
	\[ n_{(2,1,1),3} = C_{(2,1,1),(2)}= 1. \]
	Finally, $T_{p,1}=\emptyset$, so $n_{p,2}=0$. The coefficient of $\frac{1}{2}x_{(2,1,1)}$ in the expansion of the pre-Lie logarithm will thus be $\frac{22}{5}-\frac{17}{4}+\frac{1}{3}=\frac{29}{60}$.
	
\begin{remark} Recall that a surjective map from the boxes of $p$ to $[k]$
is called a \emph{semi-standard tableau} if it is weakly increasing along the rows and strictly increasing along the columns, and a \emph{standard tableau} if moreover it is bijective. A semi-standard tableau is always in $T_{p,k}$, and in fact it satisfies the stronger requirement that all boxes labeled by $h$ are contained in the first $h$ rows of $p$. As the previous example illustrates, the inclusion of semi-standard tableaux inside $T_{p,k}$ is generally strict. However, it is also clear that for $k=|p|$ we have $T_{p,|p|}=T^{st}(p)$, where the right hand side denotes the set of standard tableaux of shape $p$. 
\end{remark}	
	Let us conclude with an observation: as we just remarked, the leading coefficient $n_{p,|p|+1}$ is associated to standard tableaux of shape $p$. 
	It is also associated to the iterated brace products of maximal length in the expansion of the pre-Lie logarithm, that is the terms
	$\frac{(-1)^{n-1}}{n}f\partial(\dots(f\partial(f))\dots))$. These are also the only terms appearing (up to a scalar coefficient) in the expansion of the pre-Lie exponential. Recalling Equation (\ref{prelieexp3})
	$$\exp_\triangleleft(x_{(0)})= \sum_{n\ge0} \frac{1}{(n+1)!}\sum_{p\in\mathcal{P}_n} N_p x_p,$$
	we get:
	\begin{prop}\label{alterplexp}
		We have for the coefficients of the pre-Lie exponential: $$N_p =\frac{n_{p,|p|+1}}{p!} =\frac{1}{p!} \sum_{s\in T^{st}(p)} \prod_{h=1}^{k} C_{s_{\le h},s_{\le h-1}}.$$
	\end{prop}
	Moreover, as, given a maximal chain of partitions $(0)<_ap_1<_a\dots <_ap_{|p|}$ the corresponding sequence of tableaux is obtained by adding at each step a single cell to the previous tableau, the computation of the $C_{p_{i+1},p_i}$ is particularly simple. More precisely, write a partition $p$ in the form $p=1^{n_1}2^{n_2}\cdots k^{n_k}\cdots$, meaning that $p$ contains $k$ with multiplicity $n_k$ for all $k>0$ (in particular, $n_k=0$ for all but a finite number of $k$). Then, given $p'<p$ with $|p'|=|p|-1$ either $p'=1^{n_1-1}2^{n_2}\cdots k^{n_k}\cdots$ or there is a unique $j\ge2$ such that $p'=1^{n_1}\cdots(j-1)^{n_{j-1}+1}j^{n_j-1}\cdots k^{n_k}\cdots$: it is easy to check that in the first situation $C_{p,p'}=|p|-l(p)+1$ and in the second situation $C_{p,p'}=j(n_{j-1}+1)$. This provides another way to implement the computation of the exponential coefficient $N_p$.
	
	\section{The pre-Lie flow map} \label{sec:flow} 
	
	Let us conclude with a study of the pre-Lie flow map. We dot not discuss convergence issues and treat expressions formally. One may for example assume that the pre-Lie algebras considered below are graded connected or, more generally, equipped with a decreasing filtration $\dots L_i\subset L_{i-1}\subset \dots L_1=L$ compatible with the pre-Lie product and such that $\bigcap_iL_i=0$. 
	
	Given such a pre-Lie algebra $(L,\curvearrowleft)$, the set of group-like elements in its enveloping algebra $(\Q[L],\ast)$ is a group (a general property of Hopf algebras), and one can transport this group law to $L$ by the $\log^\odot$ map. For $v,w\in L$ we get, denoting $\circledcirc$ the product:
	$$v\circledcirc w:=\log^\odot(\exp^\odot(v)\ast \exp^\odot (w)).$$
	One can show that, equivalently,
	$$v\circledcirc w=w+v\{\exp^\odot(w)\}.$$
	This allows to associate to $v$ the group $G(v):=\{v^{\circledcirc n},n\in\mathbb{Z}\}$, with $v^{\circledcirc 1}=v$,
	$$v^{\circledcirc n}=v+v^{\circledcirc n-1}\{\exp^\odot (v)\}\qquad\mbox{for $n>1$},$$
	which unravels, after an easy induction, to
	\begin{equation}\label{flowfla}
		v^{\circledcirc n}=nv+\sum_{i=1}^{n-1}\binom{n}{i+1}\sum_{k_1,\ldots,k_i\ge1}\frac{1}{k_1!\cdots k_i!} (\dots ((v\{v^{\odot k_1}\})\dots )\{v^{\odot k_i}\},
	\end{equation}
	where we recognize, up to scalar coefficients, terms of the expansion of the pre-Lie logarithm.
	
	This group embeds into the one-parameter group associated to the flow $F(v):=\{v^{\circledcirc t},t\in\R\}$, where, as  $\exp_\curvearrowleft=\log^\circ\circ\exp^\ast$ and $\log_\curvearrowleft=\log^\ast\circ\exp^\odot$
	$$v^{\circledcirc t}:=\exp_\curvearrowleft(t\log_\curvearrowleft (v)).$$
	
	We will focus on the generic case for Novikov algebras, namely the ``formal flow group'' $F(x)$, where $x$ stand as usual for the generator of the free Novikov algebra. We write
	$$Q(t):=x_{(0)}^{\circledcirc t}:=\exp_\tl (t\log_\tl (x_{(0)})).$$
	Notice that we can recover the pre-Lie logarithm $\log_\triangleleft(x_{(0)})$ from the formal series $Q(t)$ according to 
	\[ \log_\triangleleft(x_{(0)})= Q'(0). \]
	
	Given a partition $p=(p_1,\ldots,p_k)$, we expand $Q(t)$ similarly to the expansion of the pre-Lie logarithm, that is according to 
	\[ Q(t) = \sum_p Q_p(t) \frac{x_p}{p!}\]
	
	The link between the expansion (\ref{flowfla}) of the $x_{(0)}^{\circledcirc n}$ leads then to an expansion of $Q(t)$ in the binomial basis $\binom{t}{k}$:
	\begin{prop}\label{flowbasicfla} We have, for the expansion of the formal flow $Q(t)$:
		\[ Q_p(t) = \sum_{k=1}^{|p|+1} n_{p,k}\binom{t}{k}. \]
	\end{prop}
	\begin{proof} Denote by $q(t)$ the polynomial on the right hand side of the desired identity.
		
	With the notations from the previous section (and in particular abbreviating $x_{(0)}$ with $x$), recalling the definition of $n_{p,i+1}$ as the coefficient of $\frac{x_p}{p!}$ in the expansion of \eqref{eq:npi}, for all $n\ge1$ we have
		\begin{eqnarray} \nonumber \sum_{p}Q_p(n)\frac{x_p}{p!} = x^{\circledcirc n} &=& nx+\sum_{i=1}^{n-1}\binom{n}{i+1}\sum_{k_1,\ldots,k_i\ge1}\frac{1}{k_1!\cdots k_i!} (\cdots(x\{x^{\odot k_1}\})\cdots)\{x^{\odot k_i}\}\\\nonumber
			&=& \sum_{p}\left(\sum_{i=0}^{n-1}\binom{n}{i+1}n_{p,i+1}\right)\frac{x_p}{p!}
		\end{eqnarray}
	and thus	
		\begin{equation}
			Q_p(n) = \sum_{i=0}^{n-1} \binom{n}{i+1}n_{p,i+1} =q(n).
		\end{equation}
	Since the above identity was proved for all $n>0$, this implies $Q_p(t)=q(t)$ as desired.	\end{proof}
	As a concrete example, for the partition $p=(2,2)$ the tableaux in $T_{p,4}, T_{p,3}, T_{p,2}$ are respectively
	\begin{multline*} T_{p,4}=\left\{\begin{array}{cc} 1 & 2 \\ 3 & 4
		\end{array}\:,\quad\begin{array}{cc} 1 & 3 \\ 2 & 4
		\end{array}\right\}\\ T_{p,3}=\left\{\begin{array}{cc} 1 & 1 \\ 2 & 3
		\end{array}\:,\quad\begin{array}{cc} 1 & 2 \\ 2 & 3
		\end{array}\:,\quad \begin{array}{cc} 1 & 2 \\ 3 & 3
		\end{array}\:,\quad\begin{array}{cc} 1 & 3 \\ 2 & 3
		\end{array}\right\}\\ T_{p,2}=\left\{\begin{array}{cc} 1 & 1 \\ 2 & 2
		\end{array}\:,\quad\begin{array}{cc} 1 & 2 \\ 2 & 2
		\end{array}\right\}
	\end{multline*}
	whereas $T_{p,1}=\emptyset$, which implies $n_{p,2}=0$,
	\[ n_{(2,2),5} = C_{(2,2),(2,1)}C_{(2,1),(2)}C_{(2),(1)}+C_{(2,2),(2,1)}C_{(2,1),(1,1)}C_{(1,1),(1)} = 2\cdot 2\cdot2+2\cdot4\cdot1=16,\]
	\begin{eqnarray}\nonumber n_{(2,2),4}&=&C_{(2,2),(2,1)}C_{(2,1),(2)}+C_{(2,2),(2,1)}C_{(2,1),(1)}+C_{(2,2),(2)}C_{(2),(1)}+C_{(2,2),(1,1)}C_{(1,1),(1)} \\\nonumber &=& 2\cdot 2+2\cdot 3+2\cdot2+4\cdot 1 = 18,
	\end{eqnarray}
	\[ n_{(2,2),3}= C_{(2,2),(2)}+C_{(2,2),(1)} = 2+2 = 4,\]
	and we finally get
	$$Q_{(2,2)}(t)= 16\binom{t}{5}+18\binom{t}{4}+4\binom{t}{3}.$$
	
	\ \ 
	
	We propose below a list of other properties of formal flows.
	
	\ \
	
	Let us start with a combinatorial interpretation of the integer $Q_p(n)$.
	\begin{lemma}\label{lem:combQ} The integer $Q_p(n)$ counts pairs $(T,\ell)$, where $T$ is a planar rooted tree with associated partition $\Psi(T)=p$, and $\ell:V(T)\to[n]$ is a strictly monotone map (where $V(T)$ is seen as a poset in the usual way, with the root as its minimum). 
	\end{lemma}
	\begin{proof} The Lemma follows from our developments on the pre-Lie logarithm, but can be recovered independently as follows.
		It is a well known fact that the flow $\Omega^\circ(t):=\exp_\curvearrowleft(t\log_\curvearrowleft(\bullet))$ in the free pre-Lie algebra $(p\mathcal{L}(x),\curvearrowleft)$ of rooted trees can be expanded in the form 
		\[ \Omega^\circ(t) = \sum_{T} \Omega^\circ_T(t) \frac{T}{\sigma(T)},\]
		where the sum is over (non-planar) rooted trees, seen as posets in the usual way, $\Omega^\circ_T(t)$ is the strict order polynomial counting strictly monotone maps $V(T)\to[1,t]$, and $\sigma(T)$ is the symmetry factor counting the number of automorphisms of $T$ as a (non-planar) rooted tree. The natural projection $\Psi:p\mathcal{L}(x)\to\mathcal{N}(x):T\to x_{p(T)}$ sends the flow $\Omega^\circ(t)$ in $p\mathcal{L}(x)[[t]]$ to the flow $Q(t)$ in $\mathcal{N}(x)[[t]]$, hence we see that
		\[ Q_p(t) = \sum_{T\mbox{ s.t. }\Psi(T)=p} \frac{p!}{\sigma(T)} \Omega^\circ_T(t). \]
		Finally, the Lemma follows if we show that $\frac{\Psi(T)!}{\sigma(T)}$ counts the number $\operatorname{Pl}(T)$ of planar embeddings of $T$. This can be seen via a straightforward induction on the number of vertices, observing that for a tree $T$ of the form 
		\[ T = 
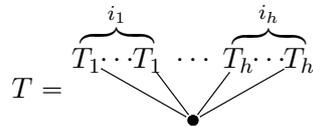
\xy {\ar@{-}(0,-4)*{\bullet};(-14,4)*{T_1}};
		{\ar@{-}(0,-4)*{\bullet};(-6,4)*{T_1}};
		{\ar@{-}(0,-4)*{\bullet};(14,4)*{T_h}};
		{\ar@{-}(0,-4)*{\bullet};(6,4)*{T_h}};
		{(0,4)*{\cdots}};{(-10,4)*{\cdots}};{(10,4)*{\cdots}};
		{(-10,8.8)*{\overbrace{\,\,\,\,}^{i_1}}};
		{(10,8.8)*{\overbrace{\,\,\,\,}^{i_h}}} \endxy \]
		(with $T_i\neq T_j$ if $i\neq j$) we have 
		\[ \operatorname{Pl}(T) = \frac{(i_1+\cdots+i_h)!}{i_1!\cdots i_h!}\operatorname{Pl}(T_1)^{i_1}\cdots \operatorname{Pl}(T_h)^{i_h}, \]	
		\[ \Psi(T)! = (i_1+\cdots+i_h)!\Psi(T_1)!^{i_1}\cdots \Psi(T_h)!^{i_h}, \]
		\[ \sigma(T) = i_1!\cdots i_h!\sigma(T_1)^{i_1}\cdots \sigma(T_h)^{i_h}. \]	
		
	\end{proof}
	
	A second remark is that $Q(t)$ satisfies two difference equations, that can be used to compute it recursively, and thus also to compute the coefficients of the pre-Lie logarithm, alternatively to the enumeration formulas we introduced in the previous section.
	Indeed, we have 
	\[Q(t+s) =Q(t)\circledcirc Q(s)\] for all $s,t\in\R$. Moreover $Q(1) = \exp_\triangleleft(\log_\triangleleft(x_{(0)}))= x_{(0)}$. In particular 
	\begin{multline*} Q(t +1 ) = Q(1))\circledcirc Q(t)=\\=x_{-1} + Q(t) + \partial(x_{-1})Q(t)+\cdots+\frac{1}{n!}\partial^n(x_{-1})Q(t)^n +\cdots =\\= Q(t) + \sum_{n\ge0}\frac{1}{n!} x_{n-1}Q(t)^n,
	\end{multline*}
	leading to the difference equation 
	\begin{equation}\label{recursion1forlog}
		\Delta Q(t) = \sum_{n\ge0}\frac{1}{n!} x_{n-1}Q(t)^n.   
	\end{equation}
	Recall that the formal inverse of the difference operator $\Delta$ acting on polynomials via
	$\Delta p(t):=p(t+1)-p(t)$
	is the indefinite sum operator
	\[\sum_{\tau=0}^{t-1}:\K[\tau] \to \K[t]: \tau^n\to\sum_{\tau=0}^{t-1} \tau^n = \frac{B_{n+1}(t)-B_{n+1}}{n+1}, \]
	where $B_n(t),B_n$ are the Bernoulli polynomials and numbers respectively. Thus we can solve \eqref{recursion1forlog} recursively up to any desired order. For instance
	\[ \Delta Q(t) = x_{-1} + \mathcal{N}(x)^{\ge2} \]
	implies $Q(t) = tx_{-1} + \mathcal{N}(x)^{\ge2}$. Substituting into \eqref{recursion1forlog} we get
	\[\Delta Q(t) = x_{-1} + x_0\Big(tx_{-1}+\mathcal{N}(x)^{\ge2}\Big)+ \mathcal{N}(x)^{\ge3} = x_{-1} + tx_0x_{-1} + \mathcal{N}(x)^{\ge3} \] and thus $Q(t) = tx_{-1} + \binom{t}{2}x_0x_{-1}+\mathcal{N}(x)^{\ge3}$. Continuing like this
	\begin{multline*} \Delta Q(t) = x_{-1} + x_0\left(tx_{-1}+\binom{t}{2} x_0x_{-1} + \mathcal{N}(x)^{\ge3}\right)+\frac{1}{2}x_1\Big(tx_{-1}+\mathcal{N}(x)^{\ge2}\Big)^2+ \mathcal{N}(x)^{\ge4} = \\ = x_{-1} + tx_0x_{-1}+ \binom{t}{2} x_0^2x_{-1} + \frac{t^2}{2} x_0x_{-1}^2+\mathcal{N}(x)^{\ge4}
	\end{multline*}
	and thus $Q(t)=tx_{-1}+\binom{t}{2}x_0x_{-1} + \binom{t}{3} x_0^2x_{-1} +\frac{B_3(t)}{6}x_1x_{-1}^2+\mathcal{N}(x)^{\ge4}$. In this way, we recover the first few term of the pre-Lie logarithm
	\[ \log_\triangleleft(x_{(0)}) = Q'(0) =  x_{(0)}-\frac{1}{2}x_{(1)}+\frac{1}{3}x_{(1,1)}+\frac{1}{12}x_{(2)}+\mathcal{N}(x)^{\ge4}.\] 
	
	More generally, from \eqref{recursion1forlog} we get the following recursion for the polynomials $Q_p(t)$: 
	\begin{prop} We have, for the formal flow:
		\begin{equation}\label{eq:rec1dual} \Delta Q_p(t) = \sum_{\stackrel{x_p = x_{k-1}x_{p_1}^{i_1}\cdots x_{p_j}^{i_j}}{i_1+\cdots+i_h=k}} \frac{k!}{i_1!\cdots i_h!} Q_{p_1}(t)^{i_1}\cdots Q_{p_h}(t)^{i_h},
		\end{equation}
		where the sums runs over (unordered) factorizations $x_p = x_{k-1}x_{p_1}^{i_1}\cdots x_{p_h}^{i_h}$ with $k,i_1,\ldots,i_h\ge1$, $i_1+\cdots+i_h=k$. 
	\end{prop}
	
	For instance, for the partition $p=(3,2,1)$ we have
	\begin{eqnarray} \nonumber x_{(3,2,1)} = x_2x_1x_0x_{-1}^4 &=& x_0(x_2x_1x_{-1}^4) = x_0x_{(3,2)}
		\\ \nonumber &=& x_1(x_2x_0x_{-1}^3)(x_{-1}) = x_1x_{(3,1)}x_{(0)} \\ \nonumber &=& x_1(x_2x_{-1}^3)(x_0x_{-1}) = x_1x_{(3)}x_{(1)} \\ \nonumber &=& x_2(x_1x_0x_{-1}^2)(x_{-1})^2 = x_2x_{(2,1)}x_{(0)}^2 \\ \nonumber &=& x_2(x_1x_{-1}^2)(x_0x_{-1})(x_{-1}) = x_2x_{(2)}x_{(1)}x_{(0)}\end{eqnarray} 
	and thus
	\[ \Delta Q_{(3,2,1)} = Q_{(3,2)}+2Q_{(3,1)}Q_{(0)} + 2Q_{(3)}Q_{(1)}+3Q_{(2,1)}Q_{(0)}^2+6Q_{(2)}Q_{(1)}Q_{(0)}. \]
	
	The second natural recursion  brings us back instead to the coefficients involved in the computation of the pre-Lie logarithm. It is obtained from $$Q(t+1)=Q(t)\circledcirc Q(1)=Q(t)\circledcirc x_{-1},$$ which, when interpreted using the construction of the coefficients $C_{p,p'}$ leads to 
	\begin{prop}We have, for the formal flow
		\begin{equation}\label{recursion2forlog}
			\Delta Q(t) = x_{-1}+\sum_{n\ge1}\frac{1}{n!} \partial^n\big(Q(t)\big)x_{-1}^n.   \end{equation}
		that is, for the coefficients $Q_{p}(t)$,
		\begin{equation}\label{eq:rec2} \Delta Q_{p}(t) = \sum_{p'<_a p} C_{p,p'} Q_{p'}(t),
		\end{equation}
	\end{prop}
	
	\section{Appendix} For future reference, in this Appendix we write down explicitly the recursions for the polynomials $Q_p(t)$ coming from \eqref{recursion1forlog} and \eqref{recursion2forlog} respectively, for all $p$ with $|p|\le 5$. The first identity comes from \eqref{recursion1forlog}, and the second from \eqref{recursion2forlog}.
	\[\Delta Q_{(0)} = 1\]
	\[ \Delta Q_{(1)}= Q_{(0)}\]\[ \Delta Q_{(2)} =Q_{(0)}^2 = 2Q_{(1)}+Q_{(0)}\]
	\[\Delta Q_{(1,1)} = Q_{(1)}\]
	\[\Delta Q_{(3)} =Q_{(0)}^3= 3Q_{(2)}+3Q_{(1)}+Q_{(0)}\]
	\[\Delta Q_{(2,1)} =Q_{(2)}+2Q_{(1)}Q_{(0)}= 2Q_{(2)} + 4 Q_{(1,1)}+3Q_{(1)}\]
	\[\Delta Q_{(1,1,1)} = Q_{(1,1)}\]
	\[\Delta Q_{(4)}=Q_{(0)}^4 = 4Q_{(3)}+6Q_{(2)}+4Q_{(1)}+Q_{(0)}\]
	\[ \Delta Q_{(3,1)} = Q_{(3)}+3Q_{(1)}Q_{(0)}^2=3Q_{(3)} +3Q_{(2,1)}+6Q_{(2)}+ 6Q_{(1,1)}+ 4Q_{(1)} \]\[ \Delta Q_{(2,2)}=2Q_{(2)}Q_{(0)} = 2Q_{(2,1)}+2Q_{(2)} + 4Q_{(1,1)} + 2Q_{(1)}\]
	\[\Delta Q_{(2,1,1)}=Q_{(2,1)}+2Q_{(1,1)}Q_{(0)}+Q_{(1)}^2 = 2Q_{(2,1)}+ 6Q_{(1,1,1)}+ Q_{(2)}+ 5Q_{(1,1)} \]
	\[\Delta Q_{(1,1,1,1)} = Q_{(1,1,1)}\]
	\[ \Delta Q_{(5)} =Q_{(0)}^5= 5Q_{(4)}+10Q_{(3)}+10Q_{(2)} + 5Q_{(1)}+Q_{(0)}\] \[ \Delta Q_{(4,1)}=Q_{(4)}+4Q_{(1)}Q_{(0)}^3 =4Q_{(4)} + 4Q_{(3,1)}+12Q_{(3)} + 6Q_{(2,1)}+12Q_{(2)}+8Q_{(1,1)}+5Q_{(1)}\]
	\[ \Delta Q_{(3,2)}=2Q_{(3)}Q_{(0)}+3Q_{(2)}Q_{(0)}^2 = 2Q_{(3,1)} +6Q_{(2,2)}+3Q_{(3)} +  9Q_{(2,1)}+8Q_{(2)}+12Q_{(1,1)}+5Q_{(1)} \]
	\begin{equation*} \Delta Q_{(3,1,1)} = Q_{(3,1)}+3Q_{(1,1)}Q_{(0)}^2+3Q_{(1)}^2Q_{(0)}= 3Q_{(3,1)} + 3Q_{(2,1,1)}+3Q_{(3)}+6Q_{(2,1)} + 9Q_{(1,1,1)}+3Q_{(2)}+7Q_{(1,1)}
	\end{equation*}\[ \Delta Q_{(2,2,1)}=Q_{(2,2)}+2Q_{(2,1)}Q_{(0)}+2Q_{(2)}Q_{(1)} = 3Q_{(2,2)} + 4Q_{(2,1,1)}+6Q_{(2,1)} + 12Q_{(1,1,1)}+2Q_{(2)}+8Q_{(1,1)}\]
	\[ \Delta Q_{(2,1,1,1)}=Q_{(2,1,1)}+2Q_{(1,1,1)}Q_{(0)}+2Q_{(1,1)}Q_{(1)} = 2Q_{(2,1,1)} + 8Q_{(1,1,1,1)}+Q_{(2,1)} + 7Q_{(1,1,1)}\]
	\[ \Delta Q_{(1,1,1,1,1)} =Q_{(1,1,1,1)}\]
	We also write down the expansion of $Q_p(t)$ (for $|p|\le5$) in the bases $t^n,\binom{t}{n}$ of $\Q[t]$. 
	\[ Q_{(0)}=t = \binom{t}{1}\]
	\[ Q_{(1)} = \frac{1}{2}t^2-\frac{1}{2}t = \binom{t}{2}\]
	\[ Q_{(2)} = \frac{1}{3}t^3-\frac{1}{2}t^2+\frac{1}{6}t = 2\binom{t}{3}+\binom{t}{2}\]
	\[ Q_{(1,1)} = \frac{1}{6}t^3-\frac{1}{2}t^2+\frac{1}{3}t=\binom{t}{3}\]
	\[Q_{(3)}= \frac{1}{4}t^4-\frac{1}{2}t^3+\frac{1}{4}t^2 = 6\binom{t}{4}+6\binom{t}{3}+\binom{t}{2}\]
	\[ Q_{(2,1)} = \frac{1}{3} t^4-\frac{7}{6}t^3+\frac{7}{6}t^2-\frac{1}{3}t = 8\binom{t}{4}+5\binom{t}{3}\]
	\[ Q_{(1,1,1)} = \frac{1}{24}t^4-\frac{1}{4}t^3+\frac{11}{24}t^2-\frac{1}{4}t = \binom{t}{4}\]
	\[ Q_{(4)} = \frac{1}{5}t^5-\frac{1}{2}t^4+\frac{1}{3}t^3-\frac{1}{30}t=24\binom{t}{5}+36\binom{t}{4}+14\binom{t}{3}+\binom{t}{2}\]
	\[ Q_{(3,1)} = \frac{7}{20}t^5-\frac{11}{8}t^4+\frac{5}{3}t^3-\frac{5}{8}t^2-\frac{1}{60}t =42\binom{t}{5}+51\binom{t}{4}+13\binom{t}{3}\]
	\[ Q_{(2,2)} = \frac{2}{15} t^5 -\frac{7}{12}t^4+\frac{5}{6}t^3-\frac{5}{12}t^2+\frac{1}{30}t= 16\binom{t}{5}+18\binom{t}{4}+4\binom{t}{3} \]
	\[ Q_{(2,1,1)} = \frac{11}{60}t^5-\frac{9}{8}t^4+ \frac{7}{3}t^3 -\frac{15}{8}t^2+\frac{29}{60}t=  22\binom{t}{5}+17\binom{t}{4}+\binom{t}{3} \]
	\[ Q_{(1,1,1,1)} = \frac{1}{120}t^5-\frac{1}{12}t^4+\frac{7}{24}t^3-\frac{5}{12}t^2+\frac{1}{5}t =\binom{t}{5} \]
	\[ Q_{(5)} = \frac{1}{6}t^6-\frac{1}{2}t^5+\frac{5}{12}t^4-\frac{1}{12}t^2 =  120\binom{t}{6}+240\binom{t}{5}+150\binom{t}{4}+30\binom{t}{3}+\binom{t}{2}\]
	\[ Q_{(4,1)} =\frac{11}{30} t^6 - \frac{8}{5}t^5 + \frac{9}{4}t^4-t^3-\frac{7}{60}t^2 + \frac{1}{10}t = 264\binom{t}{6}+468\binom{t}{5}+242\binom{t}{4}+33\binom{t}{3}\]\[ Q_{(3,2)} = \frac{1}{4}t^6 - \frac{5}{4}t^5 + \frac{17}{8}t^4-\frac{4}{3}t^3 +\frac{1}{8}t^2+\frac{1}{12}t = 180\binom{t}{6}+300\binom{t}{5}+141\binom{t}{4}+16\binom{t}{3}\]
	\[ Q_{(3,1,1)} = \frac{4}{15}t^6 - \frac{67}{40}t^5+\frac{89}{24}t^4-\frac{27}{8}t^3+\frac{41}{40}t^2+\frac{1}{20}t= 192\binom{t}{6}+279\binom{t}{5}+103\binom{t}{4} + 6\binom{t}{3}\]
	\[ Q_{(2,2,1)} = \frac{17}{90} t^6 -\frac{79}{60}t^5 +\frac{119}{36}t^4-\frac{43}{12}t^3+\frac{271}{180} t^2 -\frac{1}{10}t= 136\binom{t}{6}+182\binom{t}{5}+58\binom{t}{4}+2\binom{t}{3} \]
	\[Q_{(2,1,1,1)} = \frac{13}{180} t^6 -\frac{27}{40}t^5 +\frac{169}{72}t^4-\frac{89}{24}t^3+\frac{929}{360} t^2 -\frac{37}{60}t =52\binom{t}{6}+49\binom{t}{5}+7\binom{t}{4}\]
	\[ Q_{(1,1,1,1,1)} = \frac{1}{720} t^6 -\frac{1}{48}t^5 +\frac{17}{144}t^4-\frac{5}{16}t^3+\frac{137}{360} t^2 -\frac{1}{6}t = \binom{t}{6}\]

	\vspace{1cm}
	\noindent {\bf Acknowledgements.} 
	FP acknowledges support from the grant ANR-20-CE40-0007 and the ANR -- FWF project PAGCAP.
	He thanks Sapienza Universit\`a di Roma for its hospitality.
	
	\bibliographystyle{amsplain}
	\bibliography{biblio}
	
	\noindent
	Ruggero Bandiera\\
	Dipartimento di Matematica\\ Sapienza Universit\`a
	di Roma\\ P.le A. Moro 2\\ 00185, Roma, Italy\\
	
	\noindent
	Fr\'ed\'eric Patras\\
	Universit\'e C\^ote d'Azur et CNRS\\
	UMR 7351 LJAD\\
	Parc Valrose\\
	06108 Nice Cedex 02,
	France

\end{document}